\documentclass[3p, number,sort&compress,times,fleqn]{elsarticle}

\usepackage{amsmath}
\usepackage{amssymb}
\usepackage{mathtools}
\usepackage[colorlinks = true,
            linkcolor = blue,
            urlcolor  = blue,
            citecolor = blue,
            anchorcolor = blue]{hyperref}

\usepackage{float}

\usepackage{csml}

\graphicspath{{./figs/}}

\begin{document}

\begin{frontmatter}

\title{Interrogation of spline surfaces with application to isogeometric design and analysis of lattice-skin structures}

\author[engCam]{Xiao Xiao}
\author[engCam,geomLtd]{Malcolm Sabin}
\author[engCam]{Fehmi Cirak\corref{cor1}}
\ead{f.cirak@eng.cam.ac.uk}

\cortext[cor1]{Corresponding author}

\address[engCam]{Department of Engineering, University of Cambridge, Trumpington
Street, Cambridge, CB2 1PZ, UK}
\address[geomLtd]{Numerical Geometry Ltd, 19 John Amner Close, Ely, CB6 1DT, UK}

\begin{abstract}
A novel surface interrogation technique is proposed to compute the intersection of curves with spline surfaces in isogeometric analysis. The intersection points are determined in one-shot without resorting to a Newton--Raphson iteration or successive refinement. Surface-curve intersection is required in a wide range of applications, including contact, immersed boundary methods and lattice-skin structures, and requires usually the solution of a system of nonlinear equations. It is assumed that the surface  is given in form of a spline, such as a {NURBS}, {T-spline} or Catmull--Clark subdivision surface, and is convertible into a collection of B\'ezier patches. First, a hierarchical bounding volume tree is used to efficiently identify the B\'ezier patches with a convex-hull intersecting the convex-hull of a given curve segment. For ease of implementation convex-hulls are approximated with  k-dops (discrete orientation polytopes). Subsequently, the intersections of the identified B\'ezier patches with the curve segment are determined with a matrix-based implicit representation leading to the computation of a sequence of small singular value decompositions (SVDs). As an application of the developed interrogation technique the isogeometric design and analysis of lattice-skin structures is investigated. Although such structures have been common in large-scale civil engineering, current additive manufacturing, or 3d printing, technologies  make it possible to produce up to metre size lattice-skin structures with designed geometric features reaching down to submillimetre scale. The skin is a spline surface that is usually created in a computer-aided design (CAD) system and the periodic lattice to be fitted consists of unit cells, each containing a small number of struts. The lattice-skin structure is generated by projecting selected lattice nodes onto the surface after determining the intersection of unit cell edges with the surface. For mechanical analysis, the skin is modelled as a Kirchhoff--Love thin-shell and the lattice as a pin-jointed truss. The two types of structures are coupled with a standard Lagrange multiplier approach. 
\end{abstract}


\end{frontmatter}

\section{Introduction}
%

Almost all geometry models in isogeometric analysis are based on parametric rational spline curves and spline surfaces, such as NURBS, T-splines or subdivision surfaces~\cite{Hughes:2005aa, Bazilevs:2010aa, Cirak:2000aa, wei2015truncated, zhang2018subdivision}. The process of interacting with and exploring such geometries is commonly referred to as {\em shape interrogation} in computer-aided design (CAD) literature. Shape interrogation of curves and surfaces includes, for instance, the plotting of their differential geometric properties, like iso-curvature lines or geodesics, and the computing of their intersections with rays and planes as in ray tracing and contouring. Due to its importance in CAD there is copious literature and a number of ingenious techniques on shape interrogation, see e.g. the monographs~\cite{farin2002handbook, patrikalakis2009shape} and the references therein.  All shape interrogation involves usually, at some level, the computation of the intersection between curves, surfaces with curves, and surfaces with surfaces. The focus of this paper is on curve-surface intersection, but the presented techniques can be extended to curve-curve and surface-surface intersections. Evidently, interrogation of spline surfaces is critical in a number of isogeometric analysis applications as well. Robust intersection computations are essential, for instance, for enforcing non-penetration constraints between spline curves or surfaces~\cite{temizer2011contact, de2011large}, integrating over cut-cells~\cite{Ruberg:2011aa, Schillinger:2012aa,upreti2017signed},  coupling trimmed shell patches~\cite{marussig2017review, breitenberger2015analysis, guo2018variationally}, or, as in the present paper, for coupling shells with lattices. The number of finite elements in analysis is usually significantly larger than the number of patches in usual CAD models so that the efficiency of the intersection computations becomes paramount in isogeometric analysis.

The intersection between a spline surface and a curve leads to a nonlinear root-finding problem, which is easy to formulate but extremely hard to solve robustly. In particular, multiple intersections and curves that intersect the surface tangentially cause difficulties to most intersection algorithms. The prevalent techniques to determine intersections include triangulation~\cite{piegl1998geometry}, Newton--Raphson method~\cite{toth1985ray}, subdivision~\cite{houghton1985implementation, sederberg1986comparison, farin2002handbook-sabin}, implicitisation, to be used in this paper, or a combination thereof; see also the recent review~\cite{marussig2017review} with a focus on trimming and isogeometric analysis. In triangulation based techniques the surface is approximated with a sufficiently fine triangle mesh and the curve with a polygon. The intersection between the triangle mesh and the polygon is straightforward to compute and can serve as an approximation to the true intersection point. Subdivision based methods rely on the refinability and convex-hull properties of spline surfaces. The convex-hull property guarantees that a spline surface lies always within the convex-hull of its control points, and the refinability property enables to create a hierarchy of nested patches with decreasing convex-hull size. Hence, the intersection problem can be reformulated as the intersection between two bounding volume hierarchies for the convex-hulls of the spline patches and the curve. To simplify the intersection computation between convex-hulls it is expedient to approximate them with axis-aligned bounding boxes, spheres or, as in the paper, k-dops ({\em discrete orientation polytopes})~\cite{lin1998collision, klosowski1998efficient}.  In both triangulation and subdivision methods the smallest element or bounding volume size has to be chosen similar to the required accuracy for the intersection computation. Hence, the computing time and memory needed depend on the required accuracy. In contrast, in the Newton--Raphson method the patches can be relatively large and convergence is quadratic sufficiently close to the intersection points. The principal difficulties encountered with Newton--Raphson iteration, mainly divergence for not carefully selected starting points, are the same ones as in path-following in nonlinear static analysis~\cite{ramm1981strategies, wriggers2008nonlinear}.

In the present paper the intersections between spline surfaces and curves are computed with an implicitisation technique proposed by Bus\'e et al.~\cite{laurent2014implicit, shen2016line}. Prior to computing the intersections, the surface patches close to a curve are identified using a hierarchical k-dop bounding volume tree~\cite{klosowski1998efficient}. 
Subsequently, the intersections of the identified patches  with the curve are determined by computing the singular values and null vectors of an implicitisation matrix, without resorting to a Newton--Raphson iteration or subdivision refinement. In implicit form a parametric surface, \mbox{$\vec x = \vec x (\vec \theta)$} with \mbox{$\vec x \in \mathbb R^3$} and \mbox{$ \vec \theta \in \mathbb R^2$}, is represented as the zero iso-contour of a scalar function, i.e. \mbox{$h (\vec x)=0$}. Implicitisation naturally leads to algebraic geometry, which informally deals with systems of implicit polynomial equations of the form \mbox{$f_i(\vec x)=0$} with $i \in \mathbb N_{>0}$ and the geometric objects defined by them. For instance, a planar parametric curve, \mbox{$\vec x = ( x^1,\,x^2)^\trans =  (x^1 (\theta) ,\, x^2(\theta))^\trans \in \mathbb R^2$} and \mbox{$\theta \in \mathbb R^1$},  can be implicitised with resultants from algebraic geometry; see~\cite{hoffmann1989geometric, farin2002handbook-sed, cox2007ideals, sederbergCN:2012} for an  introduction to resultants. A polynomial resultant \mbox{$h(\vec x) =R(f_1(x^1),\, f_2( x^2))$} is a polynomial expression for \mbox{$f_1(x^1) = x^1 - x^1(\theta)=0$} and  \mbox{$f_2(x^2) = x^2 - x^2(\theta)=0$},   which is \mbox{$h (\vec x) = R(f_1( x^1),\, f_2(x^2))=0$}  if and only if both polynomials satisfy \mbox{$f_1(x^1)=0$} and \mbox{$f_2(x^2)=0$}; or in other terms only for points on the parametric curve. For a planar parametric curve the resultant  can be computed as  the determinant of the classical Sylvester or B\'ezout matrices. The method of moving lines introduced by Sederberg et al.~\cite{sederberg1994curve,sederberg1995implicitization} provides a related and more intuitive geometric approach to implicitisation. According to Sederberg et al. a planar parametric curve $\vec x(\theta)$ can be described as the intersection of two lines $l_1(\vec x,  \theta)=0$  and $l_2(\vec x, \theta)=0$ that move with the parameter $\theta$. The extension of the   moving lines idea to  surfaces is known to be challenging due to the presence of so-called base points and the impossibility of describing a surface with exactly three moving planes; see~\cite{ sederbergCN:2012} for difficulties arising from base points. The matrix-based implicitisation introduced by Bus\'e~\cite{laurent2014implicit} can be understood as the generalisation of the moving lines idea to surfaces. Instead with a closed-form implicit equation $h(\vec x)=0$ the surface is represented using a rectangular implicitisation matrix obtained by computing an auxiliary singular value decomposition (SVD). There is a drop in the rank of the implicitisation matrix when evaluated on the surface. This change in the rank can again be determined with SVD. In intersection computations the implicitisation matrix is used to formulate a generalised eigenvalue problem for directly computing the spatial and parametric coordinates of the intersection point. The relatively small SVD and eigenvalue problems can be solved with robust and efficient solvers available in most modern numerical linear algebra packages. For completeness, there are a few other matrix-based techniques for computing intersections~\cite{manocha1991new, manocha1997algebraic, dokken2001approximate}  which also consider SVDs and eigenvalue problems but are different from the implemented approach. For instance, in~\cite{manocha1991new, manocha1997algebraic} the spatial coordinates of the intersection point are determined by solving an eigenvalue problem which is equivalent to computing the roots of a polynomial resultant. While, in~\cite{dokken2001approximate} singular value decomposition is used to fit an implicitly defined function to a given parametric curve or surface.

\begin{figure}
\centering
	\subfloat[Lattice structure with BCC unit cells]
	{
		\includegraphics[scale=1.1]{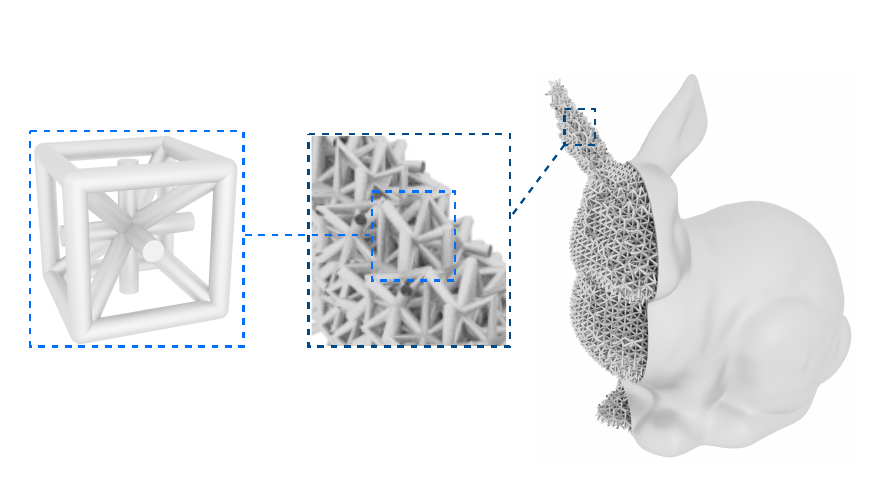}
		\label{fig:bunnyLattice}
	} 
	\hspace{0.02\linewidth}
	\subfloat[Finite element computation]
	{
		\includegraphics[scale=1.1]{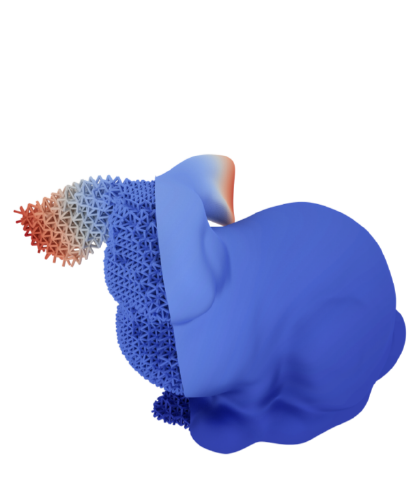}
		\label{fig:bunnyDisp}
	}
\caption{Stanford bunny with an internal lattice and its deformation under applied loading. The intersection points of the lattice with the surface are computed using the implicit matrix representation after identifying possible intersections with a k-dop bounding volume tree. }
\label{fig:bunnyLatticeSkin}
\end{figure}

As an application of the developed shape interrogation technique the isogeometric design and analysis of lattice-skin structures is considered, see Figure~\ref{fig:bunnyLatticeSkin}. The advances in additive manufacturing, or 3d printing, make it possible to produce up to metre size macroscopic lattice structures  with designed geometric features reaching down to submillimetre scale~\cite{gibson2010additive}. Structures with lattice-like components can achieve significantly higher stiffness-to-weight and strength-to-weight ratios than ordinary structures with comparable weight~\cite{gibson1999cellular, fleck2010micro}. To achieve this the struts within each unit cell must primarily be subjected to stretching (rather than bending). One of the most expedient application areas of lattice-like cellular solids is the lightweight design of structures where they are used to provide stiffness to plates and freeform shell skins~\cite{wang2013cost}. It is usually not feasible to model the entire lattice-skin structure in a CAD system due to the required very large number of geometric entities and poor scalability of CAD systems~\cite{ gibson2010additive}, although the size of lattices that can be represented is continuously increasing. In practice, this limitation is overcome by designing only the skin in the CAD system. After the skin design is completed it is triangulated and the lattice is added outside the CAD system, see e.g.~\cite{weeger2019digital}. The concurrent use of several representations leads to the usual proliferation of geometry models for representing the same part and the attendant interoperability and accuracy issues. The situation is further complicated when the generated design is to be analysed with finite elements. The interrogation technique presented in this paper allows the combination of any periodic lattice with a spline surface without the need for triangulation. After computing the intersection of the lattice with the spline surface the closest joints are projected onto the spline surface. The periodic lattice is modelled as a truss structure with pinned joints which do not transfer moments. This approximation is sufficient for lattices which are stretch dominated, see e.g.~\cite{deshpande2001collapse}. The shell skin is modelled as a Kirchhoff--Love thin-shell and is discretised with spline basis functions~\cite{Cirak:2000aa,Cirak:2011aa,zhang2018subdivision}. The truss and Kirchhoff--Love thin-shell finite elements are combined with the standard Lagrange multiplier approach, see e.g.~\cite{guidault2009bridging}. Although not pursued further here, after the design and analysis of a part is completed it can be triangulated for slicing for 3d printing. 

The outline of this paper is  as follows. In Section~\ref{sec:geometry}  the developed interrogation approach is introduced. After a cursory review of B\'ezier patches,  first  the proximity search between B\'ezier patches and curve segments using the hierarchical  k-dop bounding volume tree  is presented. Subsequently, the matrix-based implicitisation of B\'ezier patches and its use for computing their intersection points  with curve segments  are discussed.  In  Section~\ref{sec:lattice-skin},  briefly, the need and application of robust interrogation techniques in the isogeometric design and analysis of lattice-skin structures are illustrated.  The efficiency and accuracy of the developed approach are studied in~Section~\ref{sec:examples} with detailed timing studies for generating lattice-skin structure geometries and  finite element analysis of simply-supported sandwich plates. In addition, the design of an optimally stiff sandwich cap with a lattice core is studied, which demonstrates that there is an optimal material distribution between the lattice and the skin. In the Appendix, amongst others, an illustrative example for computing the intersection between two lines using the implicit matrix representation approach can be found.



%
\section{Geometry representation and interrogation \label{sec:geometry}}
%

%
\subsection{Spline surfaces \label{sec:splines}}
%
A surface is assumed to be given in form  of a spline, such as a NURBS, T-spline or Catmull--Clark subdivision surface. Each of these splines can be converted into a set of rational tensor-product B\'ezier patches, see e.g.~\cite{Piegl:1997aa, borden2011isogeometric, Stam:1998aa}. For instance, in NURBS and T-splines the surface pieces in each of the knot spans form a B\'ezier patch. In homogeneous coordinates  each B\'ezier patch is given by  
\begin{equation} \label{eq:bpatch}
	\vec{f}(\vec{\theta}) = 
	\begin{pmatrix}
		w\vec{x}\\
		w
	\end{pmatrix}
	= \sum_{\vec{i}} B^{\vec{\mu}}_{\vec{i}}(\vec{\theta})
	\begin{pmatrix}
		w_{\vec{i}}\vec{x}_{\vec{i}} \\
		w_{\vec{i}}
	\end{pmatrix} 
\end{equation}
with the coordinates
\begin{equation*}
	\vec x = (x^1, \, x^2, \, x^3)^\trans = \left ( \frac{f^1(\vec \theta)}{f^4(\vec \theta)}, \, \frac{f^2(\vec \theta)}{f^4(\vec \theta)}, \, \frac{f^3(\vec \theta)}{f^4(\vec \theta)} \right )^\trans  \in  \mathbb{R}^3 \, ,
\end{equation*}	
the bivariate B\'ezier basis functions $B^{\vec{\mu}}_{\vec{i}}(\vec{\theta}) $,  the parametric coordinates $\vec \theta = (\theta^1 , \, \theta^2) \in \mathbb{R}^2 $,  the prescribed control point coordinates $\vec x_{\vec i} = (x_{\vec{i}}^{1}, \, x_{\vec{i}}^{2}, \, x_{\vec{i}}^{3}) \in \mathbb{R}^3$ and the scalar weights $w_{\vec i} \in \mathbb R$.\footnote{Note that it may happen that \mbox{$f^1 (\vec \theta)  = f^2 (\vec \theta) = f^3 (\vec \theta) = f^4 (\vec \theta) = 0 $} for certain parameter values. These points are referred to as base points. They do not present a problem for the introduced implicitisation approach.} The B\'ezier control points $\vec x_{\vec i} $  are obtained from the spline control points by means of  a projection, referred to as B\'ezier extraction in case of tensor-product B-splines~\cite{borden2011isogeometric}.  The basis functions and the control points are labelled with the multi-index $\vec i=(i^1, \,  i^2)$. The tensor-product B\'ezier basis functions of bi-degree $\vec \mu=(\mu^1, \, \mu^2)$ are defined by
\begin{equation}
	B^{\vec{\mu}}_{\vec{i}}(\vec{\theta}) = B_{i^1}^{\mu^1}(\theta^1)B_{i^2}^{\mu^2}(\theta^2)
\end{equation}
with the Bernstein polynomials
\begin{equation} \label{eq:bernstein}
	B_i^\mu(\xi) = \binom{\mu}{i - 1}\xi^{i - 1}(1 - \xi)^{\mu - i + 1} \, , 
\end{equation}
where $  \xi \in [ 0, \, 1] $ and $ i \in \{1, \, \dotsc , \,  \mu+1 \}  $.  The Bernstein polynomials are nonnegative and sum up to one.  In applications the same degree $\mu^1 = \mu^2$ is usually used in both directions and in~\eqref{eq:bpatch} the weights are $w_{\vec i} \equiv 1$ except when quadric surfaces, like cylinders or spheres, are represented.  

Although it would be feasible to re-express a spline surface using any  polynomial basis, the Bernstein polynomials and the resulting B\'ezier patches have the advantage of numerical stability and convex-hull property. The convex-hull property guarantees that the  B\'ezier surface $\vec x(\vec \theta) $  lies always in the convex hull of its control points $\vec x_{\vec i}$ (when, as usual, $w_{\vec i} > 0$)~\cite{Farin:2002aa}. This is a direct consequence of the nonnegativity and partition of unity property of the B\'ezier basis functions.  As will be discussed in Section~\ref{sec:surfaceLineIntersect}, the intersection of  convex-hulls can be conveniently used as a preliminary check to identify possible intersections. 

It is worth noting that  in case of Catmull--Clark subdivision  the nested nature of spline rings around extraordinary vertices, with other than four connected faces, gives rise to  rings of  B\'ezier patches of decreasing size~\cite{Stam:1998aa,zhang2018subdivision}. The size of the B\'ezier patches tends to zero when the extraordinary vertex is reached.

%
\subsection{Surface-curve intersection} \label{sec:surfaceLineIntersect}
%
The computation of the intersection point between a B\'ezier patch and a given curve is an important high-level primitive task in the design of  lattice-skin structures. For both design and analysis purposes the coordinates of  the intersection point and its parametric coordinates on the surface and the curve are required. This leads to a system of nonlinear equations which are hard to solve robustly and efficiently.  Although it is possible to check the intersection of a curve with each B\'ezier patch, this has evidently a complexity of $\mathcal O(n_p)$, where $n_p$ is the number of patches,  and becomes quickly inefficient for a large number of curves.  Therefore, a hierarchical bounding volume tree and k-dops (discrete orientation polytopes) are used to first identify all the B\'ezier patches that may be intersected by the curve.  This reduces the complexity, under some mild assumptions, to  $\mathcal O(\log(n_p))$.  Subsequently,  the intersection points of the identified B\'ezier patches with the curve are  computed one by one and  valid intersections are determined.  The intersections are computed with a matrix-based implicit representation without resorting to a Newton--Raphson method.

%
\subsubsection{Bounding volume trees and k-dops for proximity search \label{sec:bvTree}}
\begin{figure}
\centering
\subfloat[Quadratic B\'ezier patch]
	{
		\includegraphics[width=0.25\linewidth]{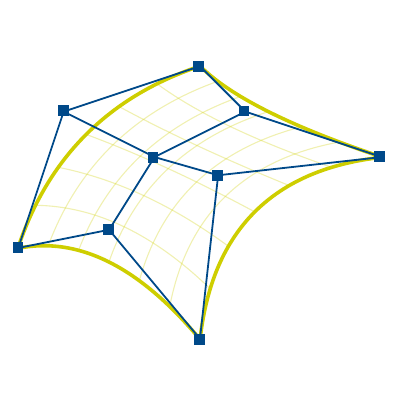}
	}
\hspace{0.025\textwidth}
\subfloat[Convex hull of the vertices]
	{
	\includegraphics[width=0.25\linewidth]{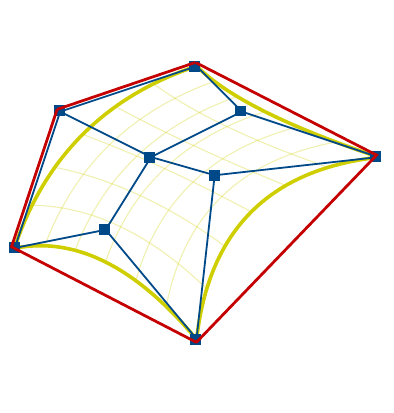}
	}
\hspace{0.025\textwidth}	
\subfloat[8-dop of the control polygon]
	{
	\includegraphics[width=0.32\linewidth]{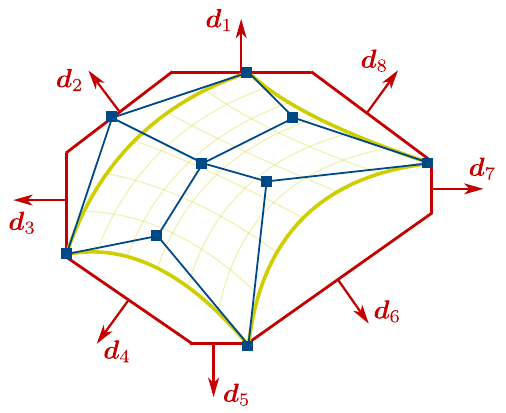}
	}
\caption{B\'ezier patches lie entirely within the convex hull and the k-dop formed by their vertices. Two-dimensional illustration with all control points $\vec x_i \in \mathbb R^2$. \label{fig:hullKdop}} 
\end{figure}

%
B\'ezier patches that may be intersected by a given curve can be efficiently identified with a hierarchical bounding volume tree.  Due to the convex-hull property of B\'ezier patches it is expedient to base  such a tree on the convex hull of patch vertices. For ease of  implementation and efficiency, however,  convex-hulls can be approximated with  k-dops~\cite{sabin1992recursive, klosowski1998efficient}, see Figure~\ref{fig:hullKdop}.  K-dops  are the generalisation of axis-aligned bounding boxes and their bounding surfaces come from a set of $k$ fixed orientations ($k \ge 6$ in 3D).  Axis-aligned bounding boxes with $k=6$ usually provide only a loose fit to convex-hulls and give many false positives. For a relatively high number of orientations  k-dops can provide a much tighter fit and their intersections are easy to compute. 

\begin{figure}
\centering
\subfloat[Root]
	{
		\includegraphics{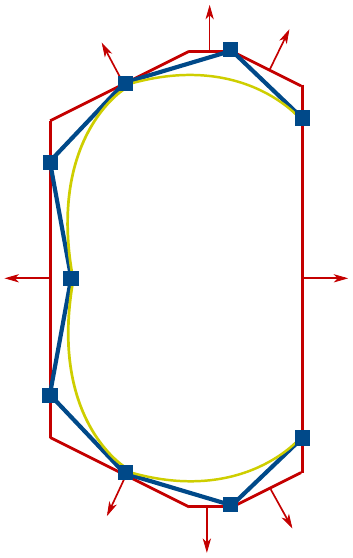}
	}
	\hspace{2cm}
\subfloat[Level 1]
	{
		\includegraphics{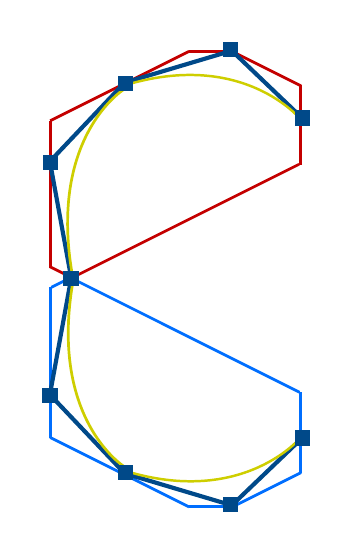}
	}
	\hspace{2cm}
\subfloat[Level 2]
	{
		\includegraphics{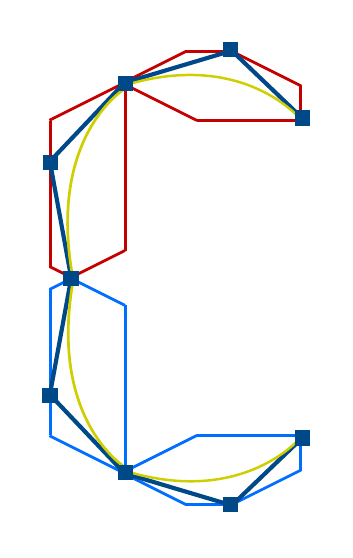}
	}
\caption{An 8-dop bounding volume tree is constructed by successively splitting the B\'ezier patches into subsets. Two-dimensional illustration for a curve consisting of four quadratic B\'ezier patches (in olive). The control polygon of each patch (in dark blue) consists of three vertices and two segments. (For interpretation of the references to colour, the reader is referred to the web version of this article.)
 \label{fig:bvTree}}
\end{figure}

In the implemented k-dop bounding volume tree each  node represents the k-dop of either one single or a set of neighbouring patches, see Figure~\ref{fig:bvTree}. The  tree is constructed by recursively splitting the patches into eight  subsets starting from the entire surface. In each step the splitting is  based on the coordinates of the patch centroids (computed as the mean of vertex coordinates). To obtain a balanced tree all the nodes on the same refinement level have to contain approximately the same number of patches. After the patches have been assigned to tree nodes, it is straightforward to determine the k-dop of a set of patches either from their individual k-dops  or  from the union of their vertices. 

\begin{figure}
\centering
\subfloat[Support heights in $\vec{d}_1$]
	{
		\includegraphics{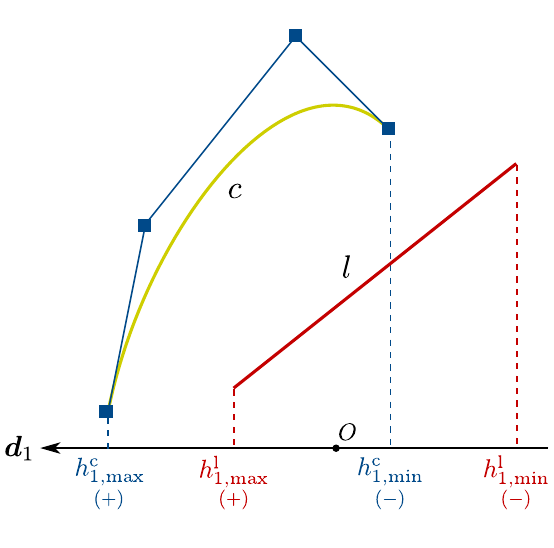}
	}
\hspace{2cm}
\subfloat[Support heights in $\vec{d}_2$]
	{
		\includegraphics{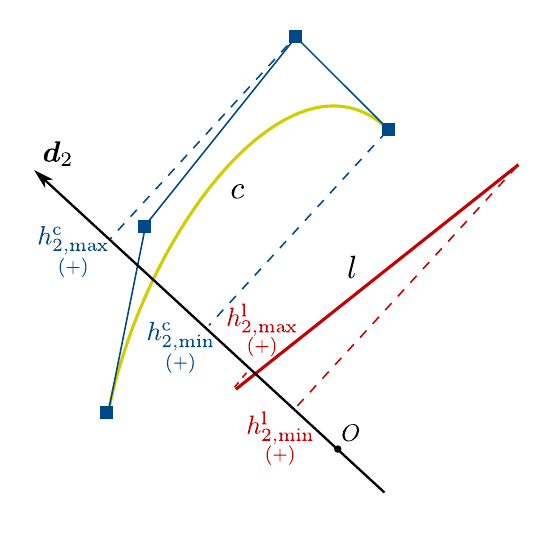}
	}
\caption{Intersection detection between a linear B\'ezier curve and a cubic B\'ezier curve. In (a) equation~\eqref{kdopsConditionForAll} is satisfied for~$j=1$. In (b) no intersection is possible because~\eqref{kdopsCondition} is satisfied.}
\label{fig:kdops}
\end{figure}

Indeed, for intersection checks there is no need to explicitly construct any k-dops.  To determine whether two  k-dops  intersect it is sufficient to consider the projections of their control vertex coordinates along the $k$ prescribed directions $\vec d_j$ with $j \in \{1, \, \dotsc , \, k\}$. As an example, in  Figure~\ref{fig:kdops} the intersection of two k-dops belonging to a linear and a cubic B\'ezier curve is illustrated. The support heights $h_{j, \text{max}}^\text{c}$ and ${h}_{j,\text{min}}^\text{c}$ of the cubic B\'ezier curve with the vertex coordinates $\vec x_i$ in the direction $\vec{d}_j$ are defined as 
\begin{equation} \label{eq:supportHeight}
h_{j, \text{max}}^\text{c} = \max_{i}(\vec{x}_{i}\cdot \vec d_j)
\quad\text{and}\quad
{h}_{j, \text{min}}^\text{c} = \min_{i}(\vec{x}_{i}\cdot{\vec d}_j) \, .
\end{equation}
The support heights $h_{j, \text{max}}^\text{l}$ and ${h}_{j,\text{min}}^\text{l}$  for the linear B\'ezier curve are determined in a similar way. The two k-dops intersect only when  
\begin{equation}  \label{kdopsConditionForAll}
\forall \, j \, ,\quad \quad  h_{j, \text{max}}^\text{c}  \geq  {h}_{j,\text{min}}^\text{l}   \quad
\text{and} \quad  h_{j, \text{max}}^\text{l}  \geq  {h}_{j,\text{min}}^\text{c}   
\end{equation}
and they do not intersect when 
\begin{equation} \label{kdopsCondition}
\exists\, j, \quad \text{such that} \quad h_{j,\max}^\text{c} < 
h_{j,\min}^\text{l}  \quad \text{or} \quad h_{j,\max}^\text{l} < 
h_{j,\min}^\text{c}  \, .
\end{equation}
Hence, in the k-dop bounding volume tree data structure it is sufficient to store only the support heights. The support heights for each tree node can be efficiently computed by starting from the leaf nodes  and determining the support heights of parent nodes from children nodes according to 
\begin{equation}
	h_{j, \text{max}}^{\text{parent}} = \max_{\text{child}} h_{j, \text{max}}^{\text {child}} 
	\quad\text{and}\quad
	h_{j, \text{min}}^{\text{parent}} = \min_{\text{child}} h_{j, \text{min}}^{\text {child}} \, .
\end{equation}

%

%
\subsubsection{Implicit matrix representation for intersection computation \label{sec:implicitisation}}
%
Next, the matrix-based implicitisation of a single B\'ezier patch $\vec f (\vec \theta)$ of bi-degree $\vec \mu = (p, \, p)$  as defined in~\eqref{eq:bpatch} is developed. To this end, consider the auxiliary vector of polynomials 
\begin{equation} \label{eq:aux}
	\vec g (\vec \theta) = \sum_{\vec j} \widetilde B_{\vec j}  (\vec \theta) \vec
	\gamma_{\vec j}
\end{equation}
with B\'ezier basis functions $\widetilde B_{\vec j} $, coefficients $\vec \gamma_{\vec j} \in \mathbb R^4$ to be yet determined and the  multi-index $\vec j = (j^1, \, j^2)$. According to~\cite{laurent2014implicit}, the bi-degree of the B\'ezier basis $\widetilde B_{\vec j}$ can be chosen with \mbox{$\widetilde{\vec \mu} = ( 2p-1, \, p-1)$} or higher; alternatively \mbox{$\widetilde{\vec \mu} = (p - 1, \, 2p - 1)$} or higher is also possible. The bi-degree of  $\widetilde B_{\vec j}$ has to be sufficiently high in order to be able to discover later all the intersection points. 
The  two vectors  $\vec f (\vec \theta)$ and $\vec g (\vec \theta)$, with each having four components, are required to be orthogonal
\begin{equation} \label{eq:orthogonality}
	\vec f (\vec \theta) \cdot \vec g (\vec \theta) = 0 \, ,
\end{equation}
that is,
\begin{equation} \label{eq:multFG}
	\left (  \sum_{\vec i} B_{\vec i} (\vec \theta ) \begin{pmatrix} w_{\vec i}
	\vec x_{\vec i} \\ w_{\vec i} \end{pmatrix} \right ) \cdot  \left ( \sum_{\vec
	j} \widetilde B_{\vec j}  (\vec \theta) \vec \gamma_{\vec j} \right ) = 0 \, .
\end{equation}
Evidently, the product of two B\'ezier basis functions $B_{\vec i} (\vec \theta)$ and $ \widetilde B_{\vec j} (\vec \theta)$ can again be expressed as a B\'ezier basis function, see Appendix~\ref{sec:productAppendix}. Hence,~\eqref{eq:multFG} can be
written with a new B\'ezier basis $\widehat B_{\vec k}$ of bi-degree~$\widehat{\vec \mu} = (3p-1, \, 2p-1)$
\begin{equation}\label{eq:implicit}
	\sum_{\vec j} \sum_{\vec k} \widehat B_{\vec k} C_{\vec k \vec j} \vec
	\gamma_{\vec j} = 0 \, ,
\end{equation} 
where $C_{\vec k \vec j}$, with the multi-indices~$\vec j$ and~$\vec k$, are the coefficients of a matrix \mbox{$\vec C \in \mathbb R^{6 p^2 \times 8p^2}$}. Note that there are  $6p^2$ basis functions  in the basis with bi-degree \mbox{$\widehat{\vec \mu} = (3p-1, \, 2p-1)$, $2p^2$} basis functions in the basis with bi-degree \mbox{$\widetilde{\vec \mu} = (2p-1, \, p-1)$} and each $\vec \gamma_{\vec j}$ has four coefficients. The matrix $\vec C$ contains the control point coordinates, associated weights and the basis transformation coefficients resulting from the B\'ezier basis multiplication.

The orthogonality constraint~\eqref{eq:orthogonality} implies $6p^2$ equations for determining $8p^2$ unknown components of $\vec{\gamma}_{\vec{j}}$'s so that the matrix $\vec C$ must be rank deficient. The non-trivial $\vec{\gamma}_{\vec{j}}$'s satisfying~\eqref{eq:implicit} must lie in the null space of $\vec C$, which is straightforward to compute with a singular value decomposition (SVD)~\cite{strangLinAlg}. The singular value decomposition of $\vec C$ reads
\begin{equation}
	\vec C = \vec U_{ c} \vec \Sigma_{c} \vec V_{c}^\trans \, 
\end{equation}
where $\vec U_{c} \in \mathbb R^{6 p^2 \times 6 p^2} $ is a matrix of left singular  vectors, $ \vec \Sigma_{c} \in \mathbb R^{6 p^2 \times 8 p^2} $ is the diagonal matrix of singular values and \mbox{$\vec V_{c} \in \mathbb R^{8 p^2 \times 8 p^2}$} is a matrix of right singular  vectors. The singular values in $ \vec \Sigma_{c} $ are usually sorted in descending order starting from the top. The right null vectors are the columns of $\vec V_c$ which correspond to zero diagonal entries in  $ \vec \Sigma_{c} $. For a specific B\'ezier patch $\vec f (\vec \theta)$ the number of right null vectors of $\vec C$ depends on its control vertex coordinates~$\vec x_{\vec i}$ and the corresponding weights~$w_{\vec i}$.

\begin{figure}
\centering
	\subfloat[scale=1.0][Six intersection lines of the four planes] 
	{
		\includegraphics{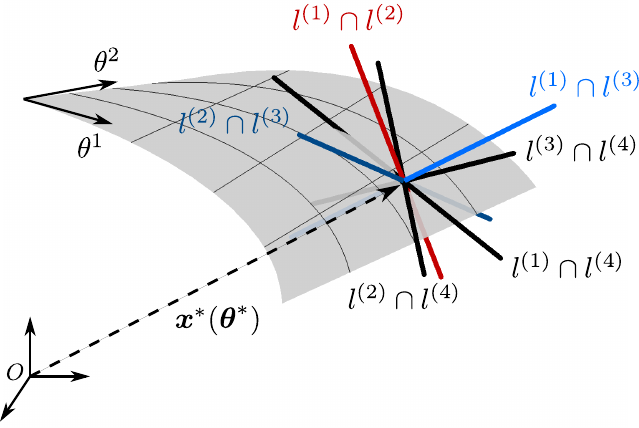}
		\label{fig:bezierSurfacePlane1}
	}
	\hspace{0.055\textwidth}
	\subfloat[scale=1.0][Three of the planes and their intersection lines]
	{
		\includegraphics{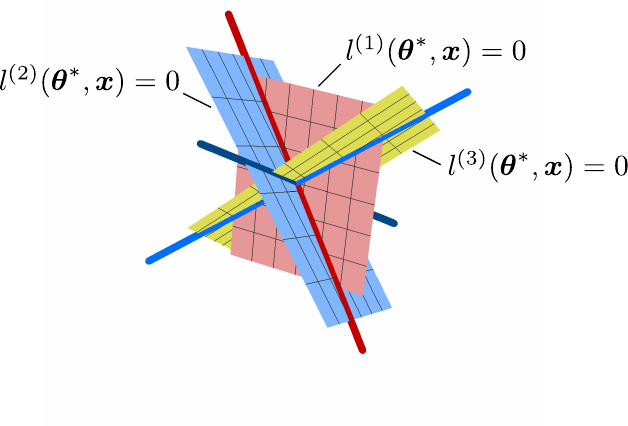}
		\label{fig:bezierSurfacePlane2}
	}
\caption{A B\'ezier patch $\vec f(\vec x)$ and its four moving planes $l^{(i)}(\vec{\theta}^*, \, \vec{x})=0$ at  $\vec x^* (\vec \theta^* ) $. For the sake of clarity, in (a) instead of the four moving planes their six pairwise intersection lines are shown and in (b) only three of the four moving planes and their three pairwise intersection lines are shown. Three of the intersection lines are present in both (a) and (b), and have in both the same colours.}
\label{fig:movingPlanes}
\end{figure}

Denoting the right null vectors of $\vec C$ with $\vec \gamma_{\vec j}^{(i)}$ and introducing them into the auxiliary vector of polynomials~\eqref{eq:aux} yields
\begin{equation} \label{eq:auxSing}
	\vec g^{(i)} (\vec \theta) = \sum_{\vec j} \widetilde B_{\vec j}  (\vec \theta)
	\vec \gamma_{\vec j}^{(i)} \, .
\end{equation}
To advance an intuitive geometric interpretation to implicitisation consider,  for a fixed parametric coordinate $\vec \theta^*$, the set of planes defined by the vectors~$\vec g^{(i)} (\vec \theta^*)$ with the implicit equations
 \begin{equation}\label{eq:planes}
 	 l^{(i)} (\vec \theta^*, \, \vec x ) = \vec g^{(i)}(\vec \theta^*) \cdot 
	\begin{pmatrix}
		w \vec x \\ w  
	\end{pmatrix}	
	= 0 \, , \quad \text{that is, } \quad
	l^{(i)} (\vec \theta^*, \, \vec x ) = 
	 \vec g^{(i)}(\vec \theta^*) \cdot 
	\begin{pmatrix}
		 \vec x \\ 1
	\end{pmatrix}  = 0 \, .
	 \end{equation}
Due to the orthogonality condition~\eqref{eq:orthogonality} the surface point with the homogeneous coordinate $\vec f(\vec \theta^*)$ and the surface coordinate
\begin{equation*}
 	\vec x^*= \left ( \frac{f^1(\vec \theta^*)}{f^4(\vec \theta^*)}, \, \frac{f^2(\vec \theta^*)}{f^4(\vec \theta^*)}, \, \frac{f^3(\vec \theta^*)}{f^4(\vec \theta^*)} \right )^\trans   
\end{equation*}
satisfies all plane equations and, hence, is present on all planes, i.e. $ l^{(i)} (\vec \theta^*,\,  \vec x^*) = 0 \, ,  \forall \, \,   i $.  For a point $\vec x^*$ to be present on two or more planes it must lie on their intersection. Note that the intersection of two planes defines a line and the intersection of a plane with a line defines a point.  These observations motivate the definition of the surface~$\vec f(\vec \theta)$ as the intersection of a set of planes.  Evidently, the planes must move, or change their inclination, while the parametric coordinate $\vec \theta$ is varied. Hence, each $ l^{(i)} (\vec \theta, \, \vec x ) =0$ for a fixed $i$ represents a family of planes and is traditionally referred to  as a {\em moving plane} or a {\em pencil of planes}~\cite{sederberg1994curve}.  Although three planes are sufficient to define a point, according to~\cite{laurent2014implicit} more than three planes must be used, in general, to describe an entire patch.  As an example in Figure~\ref{fig:movingPlanes} the description of a B\'ezier patch with four moving planes is illustrated. For the sake of clarity in~Figure~\ref{fig:bezierSurfacePlane2} only three of the planes are depicted.
  
Instead of the geometrically motivated moving plane interpretation it is possible to advance an equivalent matrix-based linear algebra explanation. To this end, the plane equations~\eqref{eq:planes} are, after introducing~\eqref{eq:aux}, rewritten as 
\begin{equation} \label{eq:planesMatrix}
 l^{(i)} (\vec \theta^*, \, \vec x )  = 
\vec{g}^{(i)}(\vec{\theta}^*)\cdot
	\begin{pmatrix}
	\vec{x} \\
	1
	\end{pmatrix} =
\sum_{\vec{j}}\widetilde{B}_{\vec{j}}(\vec{\theta}^*) \left ( \vec{\gamma}_{\vec{j}}^{(i)} \cdot
	\begin{pmatrix}
	\vec{x} \\
	1
	\end{pmatrix}
	\right )
= \sum_{\vec{j}}\widetilde{B}_{\vec{j}}(\vec{\theta}^*)M_{\vec{j}}^{(i)}(\vec{x})  = 0 \, ,
\end{equation}
where $M_{\vec{j}}^{(i)}(\vec{x})$ constitutes the $i$-th column of a matrix $\vec{M}(\vec{x})$, that is,
\begin{equation} \label{eq:MRep}
\vec{M}(\vec{x}) =
\begin{pmatrix}
M_{\vec{j}}^{(1)}(\vec{x}) & M_{\vec{j}}^{(2)}(\vec{x}) & \cdots &
\end{pmatrix} \, .
\end{equation}
The number of columns of~$\vec M(\vec x)$ is equal to the number of planes or null vectors~$\vec \gamma_{\vec j}^{(i)}$, and its number of rows is~$2 p^2$.  Again, due to the orthogonality condition~\eqref{eq:orthogonality}  for a point with the homogeneous coordinate~$\vec f(\vec \theta^*)$ and the surface coordinate~$\vec x^*$  the set of equations~\eqref{eq:planesMatrix} are satisfied for all $i$. This implies that there has to be a change in the rank of~$\vec M(\vec x)$ when evaluated at $\vec x^*$. Only at~$\vec x^*$, $\widetilde B_{\vec j} (\vec \theta^*) $ is in the null-space of~$\vec M(\vec x) $.  According to the rank-nullity theorem of linear algebra, the increase in the dimension of the null-space of~$\vec M(\vec x) $ is associated with a decrease in its rank~\cite{strangLinAlg}. Furthermore, the  minimum degree \mbox{$\widetilde {\vec \mu} = (2p-1, \, p-1)$} suggested for the basis~$\widetilde B_{\vec j}$ ensures that~$\vec M(\vec x)$ is full rank except on the surface~\cite{laurent2014implicit}.  This motivates the use of the change of rank of $\vec M(\vec x)$ to define the surface~$\vec f(\vec \theta)$, which gives rise to the notion of {\em matrix representation} of a spline surface.

The matrix $\vec M(\vec x)$ allows the formulation of the interrogation of the B\'ezier patch~$\vec f(\vec \theta)$ as a linear algebra problem.  As an example the intersection of a patch with a line $\vec r(\xi)$ is discussed next. For the intersection of a patch with a higher-order curve see Appendix~\ref{sec:linearisationAppendix}. The line is assumed to be given in parametric form with 
\begin{equation}\label{eq:ray}
	\vec r (\xi) = \vec c_1 \xi + \vec c_0 \, ,
\end{equation}
where $\vec c_1, \vec c_0 \in \mathbb R^3$ are two prescribed vectors  and $\xi \in \mathbb R$ is the parametric coordinate. It is assumed that the patch~$\vec f(\vec \theta)$ intersects the line at the parameter value(s) $\xi^*$, which can be several. To determine $\xi^*$ the  line~\eqref{eq:ray} is introduced into~\eqref{eq:MRep} yielding
\begin{equation}\label{eq:linPencil}
	\vec M( \vec r (\xi)) = \vec A - \xi \vec B  \, ,
\end{equation}
where the matrices  $\vec A$ and $\vec B$ contain the components of~$\vec M( \vec r (\xi))$ independent and linear in~$\xi$, respectively. As discussed, at the intersection points~$\vec M( \vec r (\xi))$ must be rank deficient or~\eqref{eq:planesMatrix} must be satisfied. Hence, the generalised eigenvalue problem 
\begin{equation} \label{eq:evKronecker}
	( \vec A - \xi \vec B) \vec \phi  = \vec 0 
\end{equation}
is considered for computing all the intersection points~$\xi^*$ of the line with the patch, where only the real eigenvalues are relevant. Here, the eigenvectors~$\vec \phi$ are of no significance for the intersection computations. 

The matrix~$\vec M( \vec r (\xi)) $, and hence the matrices~$\vec A$  and~$\vec B$, are usually not square. Therefore, the eigenvalue problem~\eqref{eq:evKronecker} needs to be first transformed into a block column echelon form by applying a sequence of orthogonal transformations. Note that the eigenvalues are invariant with respect to orthogonal transformations and that in block column echelon form the eigenvalues  depend only on the diagonal blocks.  Following the  Algorithm 4.1 given in~\cite{van1979computation}, first the SVD of the matrix $\vec B$ is computed
\begin{equation}
	\vec B = \vec U_b \vec \Sigma_b \vec V_b^\trans \,  .
\end{equation}
This decomposition has the usual orthogonality properties~$\vec U_b \vec U_b^\trans=\vec I$ and~$\vec V_b^\trans \vec V_b = \vec I$, where $\vec I$ is the identity matrix, and  the singular eigenvalues in~$\vec \Sigma_b$ are sorted in descending order. The matrices~$\vec B$ considered are rank deficient so that~$\vec \Sigma_b$ contains a number of zero diagonal entries.  Hence, the right singular vectors~$\vec V_b$ can be used to partition the matrices   $\vec A$  and $\vec B$ such that 
 \begin{equation} \label{eq:columnCompression}	
 	\vec B \vec V_b = 
	\vec U_b \vec \Sigma_b \vec V_b^\trans \vec V_b
	=
	\begin{pmatrix}
		\vec B_{11}' & \vec 0
	\end{pmatrix} 
	\quad 
	\text{and}
	\quad
		\vec A \vec V_b = 
	\begin{pmatrix}
		\vec A_{11}' & \vec A_{12}'
	\end{pmatrix} \, .
 \end{equation}
The number of columns of  $\vec B_{11}'$ corresponds to the rank of $\vec B$ and $\vec 0$ denotes a zero matrix.  $\vec A_{11}'$ is chosen to have the same number of columns as $\vec B_{11}'$. Next the SVD of the matrix $\vec A_{12}'$ is computed
\begin{equation}
	\vec A_{12}'= \vec U_{a_{12}} \vec \Sigma_{a_{12}} \vec V_{a_{12}}^\trans \, .
\end{equation}
The corresponding left singular vectors are used to further partition the matrices~$\vec B \vec V_{b}$ and~$\vec A \vec V_b$  such that
 \begin{equation} \label{eq:rowCompression}	
 	\vec U_{a_{12}}^\trans \vec B \vec V_{b} = 
	\begin{pmatrix}
		\vec B_{11}'' & \vec 0  \\
		\vec B_{21}'' & \vec 0
	\end{pmatrix} 
	\quad 
	\text{and}
	\quad
	\vec U_{a_{12}}^\trans	\vec A \vec V_b = 
	\begin{pmatrix}
		\vec A_{11}'' & \vec A_{12}'' \\
		\vec A_{21}'' & \vec 0
	\end{pmatrix} \, .
 \end{equation}
The number of rows of the upper blocks corresponds to the rank of $\vec A_{12}'$. With these partitioned matrices the generalised eigenvalue problem~\eqref{eq:evKronecker} can be rewritten as
 \begin{equation} \label{eq:eigenvalueTransform}	
 	\vec U_{a_{12}}^\trans ( \vec A - \xi \vec B) \vec V_b \vec \phi  = 
	\begin{pmatrix}
		\vec A_{11}''  - \xi \vec B_{11}''& \vec A_{12}'' \\
		\vec A_{21}'' -  \xi \vec B_{21}'' & \vec 0 
	\end{pmatrix} 
	\vec \phi = \vec 0  \, .
 \end{equation}
Permuting the rows yields the desired column echelon  form 
 \begin{equation} \label{eq:eigenvalueEchelon}	
	\begin{pmatrix}
		\vec A_{21}'' -  \xi \vec B_{21}'' & \vec 0  \\
		\vec A_{11}'' - \xi \vec B_{11}''& \vec A_{12}'' 
	\end{pmatrix} 
	\vec \phi = \vec 0  \, .
 \end{equation}
The eigenvalues~$\xi$ of this problem are the same as the ones for the original eigenvalue problem~\eqref{eq:evKronecker} because only orthogonal transformations have been applied. Furthermore, to find a set of eigenvalues~$\xi^*$ and eigenvectors $\vec \phi$ that satisfy~\eqref{eq:eigenvalueEchelon} it is sufficient to consider the eigenvalue problem stipulated by the top diagonal~$(\vec A_{21}'' -  \xi \vec B_{21}'') \vec \phi_1 = \vec 0$. The top diagonal may however not yet be a square matrix in which case the above sketched transformations need to be repeated. These transformations can also be realised with other methods, for example LU decomposition~\cite{thang2009curve} and QR decomposition~\cite{beelen1988improved}, which may potentially improve the computation speed.

Among the eigenvalues $\xi^*$ of~\eqref{eq:eigenvalueEchelon} only the non-complex ones correspond to intersection points with the patch. The coordinates of the intersection points~$\vec x^*$ are obtained by evaluating the line equation $\vec x^* =\vec r(\xi^*)$.  The corresponding parametric coordinates~$\vec \theta^*$ on the patch can be determined with a final singular value decomposition \mbox{$\vec M(\vec x^*)= \vec U_m \vec \Sigma_m \vec V_m^\trans$}. According to \eqref{eq:planesMatrix} the vector of the basis function values~$\widetilde B_{\vec j}( \vec \theta^*)$ must be in the left null space of~$\vec{M}(\vec{x}^*)$. Multiple left null vectors, i.e. multiple zero diagonal entries in $\vec \Sigma_m$, indicate a self intersection of the surface at~$\vec x^*= \vec r(\xi^*)$. The parametric coordinate $\vec \theta^*$ is determined from the set of nonlinear equations
\begin{equation}
	\widetilde B_{\vec j} (\vec \theta) =  N_{\vec j} \, ,
\end{equation}
where $N_{\vec j}$ represents the coefficients of one of the left null vectors.  Although these $2 p^2$ equations are nonlinear, the two components of~$\vec \theta^*$ are determined by simply considering their ratios, like 
\begin{equation}
	\frac{\widetilde B_{(j^1, \, j^2)} (\vec \theta)}{\widetilde B_{(j^1+1, \, j^2)} (\vec \theta)} =  \frac{N_{(j^1, \, j^2)}}{N_{(j^1+1, \, j^2)}} \quad \text{or} 
	\quad \frac{\widetilde B_{(j^1, \, j^2)} (\vec \theta)}{\widetilde B_{(j^1, \, j^2+1)} (\vec \theta)} =  \frac{N_{(j^1, \, j^2)}}{N_{(j^1, \, j^2+1)}} \, ,
\end{equation}
which yields two linear equations that can be trivially solved.  For instance, for a basis of degree~\mbox{$\widetilde{\vec \mu} = ( 2p-1, \, p-1)$} and the multi-index $\vec j = (1, \,  1)$ the two ratios yield the linear equations
\begin{equation}
\frac {1- \theta^1}{(2p -1)\,\theta^1} = \frac{N_{(1, \,  1)}}{N_{(2, \, 1)}}  \quad \text{and}  \quad \frac {1- \theta^2}{(p - 1)\,\theta^2}  =  \frac{N_{(1, \, 1)}}{N_{(1, \, 2)}} \, .
\end{equation}
All intersection points with the parametric coordinates within the range $[0, \, 1]\times[0, \, 1]$ are valid.

Finally, a criterion is needed to determine the rank of a matrix in computing the singular value decompositions with finite precision arithmetic.  The numerical rank~$r_{\varepsilon}$ of a matrix can be defined considering the ratios of consecutive singular values according to  
\begin{equation}
r_\varepsilon =  \min \{ k\, \vert  \,\Sigma_{k + 1,k+1}/\Sigma_{k,k} < \varepsilon, \, k = 1, \, 2, \, \dotsc \}  \, ,
\end{equation}
where $\Sigma_{k,k}$ is the $k$-th singular value (sorted in descending order) and $\varepsilon$ is a prescribed tolerance. In our computations the tolerance is chosen with~$\varepsilon=10^{-6}$.

%
\section{Lattice-skin structures \label{sec:lattice-skin}}
%

%
\subsection{Geometry generation}
\begin{figure}
\centering
		\includegraphics[width=0.6\linewidth]{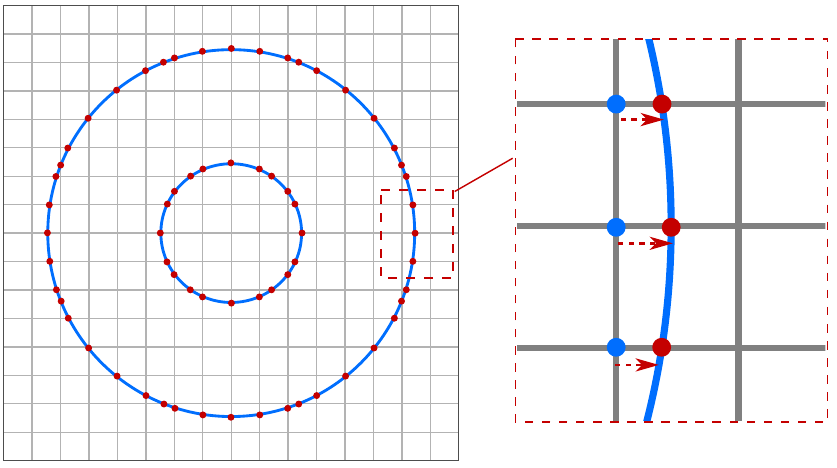}
\caption{Surface and the lattice and the projection of unit cell corners onto the surface. \label{fig:gridProject} }
\end{figure}

%
The closed surface and the specific lattice structure that is to be added to support the surface are assumed to be given, see Figure~\ref{fig:gridProject}. In this paper only periodic lattices with cubic cells of the same size are considered. Each cubic cell contains an arrangement of a small number of struts that is periodically repeatable over all cells. Depending on the placement of the joints within a cell and their connections with struts different structures are possible, like the octet, pyramidal or tetrahedral truss lattice. The assumption of cubic cells is here not overly restrictive, considering that amongst polyhedra with a small number of faces only the cube and the dodecahedron are able to uniformly fill the space~\cite{fleck2010micro}. Lattices containing several different types of polyhedral cells or with gradually changing cell size are possible but will be not considered for the sake of simplicity. 

To create the lattice-skin structure first the intersections of the lattice lines with the surface are determined. This is achieved by looping over all lines of the lattice and first determining the B\'ezier patches which might intersect the line segment and subsequently computing the exact intersection points. In design and analysis both the spatial and parametric coordinates of the intersection points are required. In most applications the size of cells is significantly smaller than the characteristic size of the surface so that there are a very large number of intersection points. 

\begin{figure}

\begin{minipage}{\textwidth}
\centering
\subfloat[Geometry fitted lattice]
	{
		\includegraphics[width=0.28\linewidth]{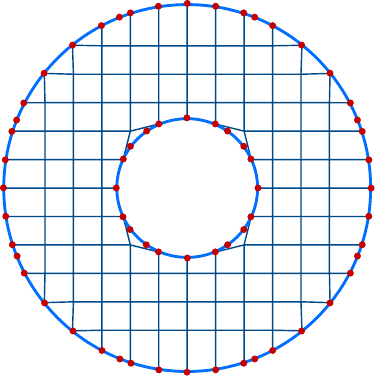}
		\label{fig:latticeFitted}
	}
\hspace{2cm}
\subfloat[Truss lattice structure]
	{
		\includegraphics[width=0.28\linewidth]{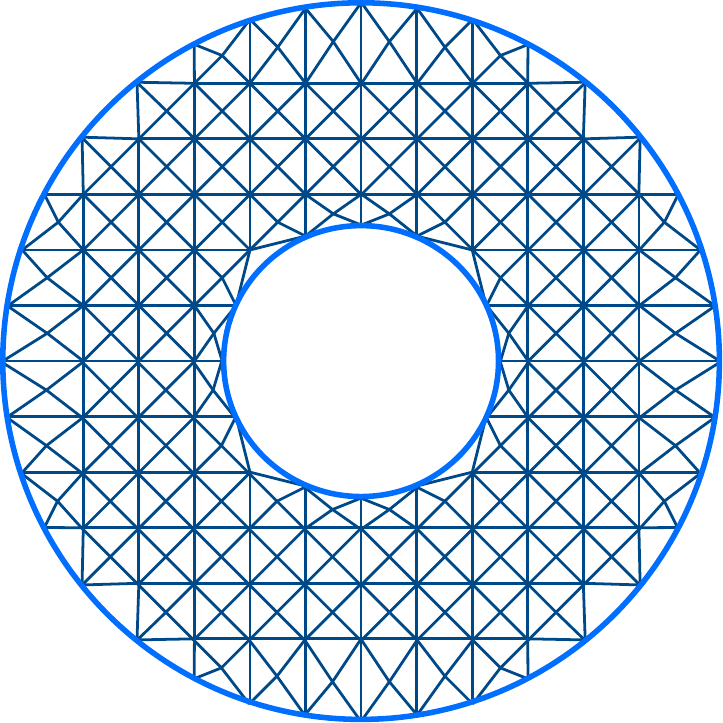}
		\label{fig:latticeGen}
	}
\end{minipage}
\caption{Truss lattice generation from a geometry fitted lattice.}
\end{figure}

After the intersection points have all been determined, the shapes of the cells intersected by the surface are modified by projecting their corner vertices onto the surface. This is achieved by first iterating over all intersection points and finding their closest lattice vertices. Subsequently, each of the identified vertices is moved to its closest intersection point by displacing it along a lattice line, see Figures~\ref{fig:gridProject} and~\ref{fig:latticeFitted}. The vertices inside the space enclosed by the surface are identified by examining their relative positions on the lattice line. Assuming that the surface is fully contained within the lattice, the vertices between the first and second, third and fourth, and so on intersection points have to be inside the surface. All the non-projected lattice vertices and their connected edges outside the closed surface are deleted.  After a lattice matching the given surface has been created, it becomes straightforward to create the truss lattice structure by first creating the joints and then the struts which connect them, see Figure~\ref{fig:latticeGen}. 

%
\subsection{Structural model}
%

%
\subsubsection{Skin}
%
We briefly review the governing equations for the linear Kirchhoff--Love model which is used for modelling the skin, for more details see, e.g.,~\cite{ciarlet2005introduction, Cirak:2000aa}. The mid-surface $\Omega^\mathrm{s}$ of the shell is parameterised with the curvilinear coordinates~$\vec \theta = (\theta^1, \, \theta^2)\in\mathbb{R}^2$. The position vector of material points on the surface is denoted with $\vec x^{\text s}(\vec \theta) \in \mathbb{R}^3$. The superscript $\text s$ has been introduced to distinguish between the variables for the shell and the lattice. 

The covariant basis vectors $\vec{a}_\alpha $ and the unit normal vector $\vec{n}$ of the mid-surface are defined with 
\begin{equation}
\vec{a}_\alpha = \frac{\partial\vec{x}^{\text s}}{\partial\theta^\alpha} =
\vec{x}^{\text s}_{,\alpha} \, , \quad
\vec{n} = \frac{\vec{a}_1\times\vec{a}_2}{|\vec{a}_1\times\vec{a}_2|} \, .
\end{equation}
Here and in the following the Greek indices take the values $\{1, 2\}$ and the summation convention is used.  As usual, the contravariant base vectors $\vec a^{\beta}$ are defined  through the relation $\vec a^{\beta} \cdot \vec a_{\alpha} = \delta^{\beta}_{\alpha}$, where  $\delta^{\beta}_{\alpha}$ is the Kronecker delta.  In case of small mid-surface displacements $\vec{u}^{\text s} (\vec \theta) \in \mathbb{R}^3 $  the  linearised membrane strain tensor $\vec{\alpha}$ and the bending strain tensor $\vec{\beta}$ are, according to~\cite{Cirak:2000aa}, derived to be 
\begin{equation}
\vec{\alpha} = \frac{1}{2}(\vec{a}_\alpha\cdot\vec{u}_{,\beta}^{\text s} +
\vec{u}_{,\alpha}^{\text s}\cdot\vec{a}_\beta) \, \vec{a}^\alpha\otimes\vec{a}^\beta
\end{equation}
and
\begin{equation}
\begin{split}
\vec{\beta} = &\Big(-\vec{u}_{,\alpha\beta}^{\text s}\cdot\vec{n} +
\frac{1}{\sqrt{a}}\left[ \vec{u}_{,1}^{\text s}\cdot(\vec{a}_{\alpha,\beta}\times\vec{a}_2) +
\vec{u}_{,2}^{\text  s}\cdot(\vec{a}_1\times\vec{a}_{\alpha,\beta}) \right] \\
&+ \frac{\vec{n}\cdot\vec{a}_{\alpha,\beta}}{\sqrt{a}}\left[
\vec{u}_{,1}^{\text s}\cdot(\vec{a}_2\times\vec{n}) +
\vec{u}_{,2}^{\text s}\cdot(\vec{n}\times\vec{a}_1)
\right]\Big) \, \vec{a}^\alpha\otimes\vec{a}^\beta \, 
\end{split} 
\end{equation}
with  $\sqrt{a} = |\vec{a}_1\times\vec{a}_2|$.

The total potential energy of the displaced thin-shell is given by 
\begin{equation} \label{eq:totPot}
\Pi^\text{s}(\vec{u}^\text{s}) = \int_{\Omega^\mathrm{s}}\!(W^m(\vec{\alpha}) + W^b(\vec{\beta}))\,\D \Omega^\mathrm{s}  - \int_{\Omega^\mathrm{s}} \vec p \cdot \vec u \, \D \Omega^\mathrm{s}\, ,
\end{equation}
where  $W^m(\vec{\alpha})$ and $W^b(\vec{\beta})$ denote the membrane and bending energy densities, and $\vec p$ is the surface load applied to the shell.  For an elastic and isotropic material the membrane and bending energy densities are given with
\begin{equation}
W^m(\vec{\alpha}) = \frac{1}{2}\frac{Et}{1 -
\nu^2}\vec{\alpha}:\vec{H}:\vec{\alpha} \, , \quad
W^b(\vec{\beta}) = \frac{1}{2}\frac{Et^3}{12(1 -
\nu^2)}\vec{\beta}:\vec{H}:\vec{\beta} \, , 
\end{equation}
where $t$ is the shell thickness, and $E$ and $\nu$ are the Young's modulus and the Poisson's ratio. The fourth-order auxiliary tensor $\vec{H}$ is defined as
\begin{equation}
 \vec H =  \left ( \nu a^{\alpha\beta}a^{\gamma\delta} +
\frac{1}{2}(1 - \nu) \left (a^{\alpha\gamma}a^{\beta\delta} +
a^{\alpha\delta}a^{\beta\gamma} \right ) \right )  \, \vec a_{\alpha} \otimes \vec a_{\beta} \otimes \vec a_{\gamma} \otimes \vec  a_{\delta}
\end{equation}
with the contravariant metric $a^{\alpha\beta} =
\vec{a}^\alpha\cdot\vec{a}^\beta$.
%

%
\subsubsection{Lattice}
%
The lattice structure is modelled  as a pin-jointed truss consisting of struts connected by pins that do not transfer moments.  Each strut deforms only by stretching and does not bend. We denote the coordinates and the displacements of the joints with $\vec x_i^{\text{l}} \in \mathbb{R}^3$ and $\vec u_i^{\text{l}} \in \mathbb{R}^3$. In case of small   displacements the axial strain in a  strut with the index $j$ is given by  
\begin{equation}
\epsilon_j = \frac{ L_{ji}  \vec{u}_{i}^\text{l} -
R_{ji} \vec{u}_{i}^\text{l}}{l_j}\cdot\vec{t}_j \, ,
\end{equation}
where $L_{ji}$ and $R_{ji}$ are two ``picking matrices'' for  selecting the displacements of the  two joints attached to the strut with the index $j$, $l_j$ is the strut length and $\vec{t}_j$ is the unit tangent to the strut. 

The total  potential energy of the lattice, with no externally applied loading at the joints, is  composed of the  potential energies of the individual struts 
\begin{equation}
\Pi^\text{l}(\vec{u}_i^{\text l}) = \sum_{j}W^a({\epsilon}_j)A_j l_j \, ,
\end{equation}
where  $A_j$ is the cross-section area of the strut with the index $j$. The external potential energy has been omitted as the loading is usually applied through the skin. For an elastic material with the Young's modulus $E$ the energy density is given by
\begin{equation}
W^a({\epsilon}_j) =
\frac{1}{2}E \epsilon_j^2 \, . 
\end{equation}
%

%
\subsubsection{Lattice-skin coupling and discretisation}
%
Some of the lattice joints $\vec x_i^{\text{l}}$ are by design located on the shell mid-surface $\vec x^\text{s}(\vec \theta)$. The displacements of these joints with~$\vec x_i^\text{l} = \vec x^{\text s} (\vec \theta_i^\text{l})$ have to be compatible with the shell displacements
\begin{equation} \label{eq:dispCompatibility}
\vec{u}^\text{l}_i = \vec{u}^\text{s}(\vec{\theta}_i^\text{l})  \, , 
\end{equation}
where $\vec{\theta}_i^\text{l}$  are the parametric coordinates of the joints  $\vec{x}_i^\text{l}$  on the surface. This condition can be imposed with Lagrange multipliers $\vec{\lambda}_i$ and considering the total  potential energy for the lattice and the skin 
\begin{equation} \label{eq:lagrangian}
	L(\vec{u}^\text{s},\, \vec{u}^\text{l},\, \vec{\lambda}) =
	\Pi^{\text{s}}(\vec{u}^\text{s}) +
	\Pi^{\text{l}}(\vec{u}^\text{l}) +
	\sum_{i}\vec{\lambda}_i\left( \vec{u}^\text{l}_i -  \vec{u}^\text{s}(\vec{\theta}_i^\text{l})  \right) \, .
\end{equation}
The  stationarity condition for this potential yields the equilibrium equations 
\begin{subequations}
	\begin{align}
	\frac{\partial\Pi^{\text{s}}(\vec{u}^\mathrm{s})}{\partial\vec{u}^\mathrm{s}}\delta\vec{u}^\mathrm{s}
	& = 0 \, , \label{eq:varshell1}\\
	\frac{\partial\Pi^{\text{s}}(\vec{u}^\mathrm{s})}{\partial\vec{u}^\mathrm{s}}\delta\vec{u}^\mathrm{s}
	- \sum_i\vec{\lambda}_i\delta\vec{u}^\mathrm{s}(\vec{\theta}_i^\mathrm{l})& = 0 \,  \quad  i \in \{\text{joints on the shell}\} \, 
	, \label{eq:varshell2}\\ 
	\frac{\partial\Pi^{\text{l}}(\vec{u}_i^\mathrm{l})}{\partial\vec{u}_i^\mathrm{l}}
	+ \vec{\lambda}_i & = 0 \, \quad  i \in \{\text{joints on the shell}\} \, ,
	\label{eq:varlattice}\\
	\frac{\partial\Pi^{\text{l}}(\vec{u}_i^\mathrm{l})}{\partial\vec{u}_i^\mathrm{l}}
	& = 0 \,  \quad  i \notin \{\text{joints on the shell}\} \, , 
	\label{eq:varlattice}\\
	\vec{u}^\text{l}_i - \vec{u}^\text{s}(\vec{\theta}_i^\mathrm{l}) & = 0 \,  \quad  i \in \{\text{joints on the shell}\} \, . 
	\label{eq:varlagrange}
	\end{align}
\end{subequations}
These equations lead to a linear system of equations after discretising the shell.  In isogeometric analysis the shell surface $\vec x^\text{s} (\theta^1, \, \theta^2)$ and displacements $\vec u^\text{s} (\theta^1, \, \theta^2) $ in~\eqref{eq:varshell1},~\eqref{eq:varshell2} and~\eqref{eq:varlagrange} are both represented with the same type of basis functions. Hence, we use the splines introduced in Section~\ref{sec:geometry} for geometry also for approximating the displacements. For displacements there is usually no need for rational splines as the exact representation of quadric sections is not  relevant. Moreover, both the geometry and displacements can be represented with B\'ezier instead of the  spline basis functions as discussed in Section~\ref{sec:splines}. If B\'ezier basis functions are used, the computed element stiffness matrices and force vectors have to be projected back to spline control points~\cite{borden2011isogeometric}. Informally, this projection ensures that the B\'ezier basis functions on neighbouring patches remain smoothly connected for the displaced structure.

%
\section{Examples \label{sec:examples}}
%

%
\subsection{Freeform surface with an internal lattice}
\begin{figure}
\begin{minipage}{\textwidth}
\centering
	\subfloat[Control mesh]
		{
			\includegraphics{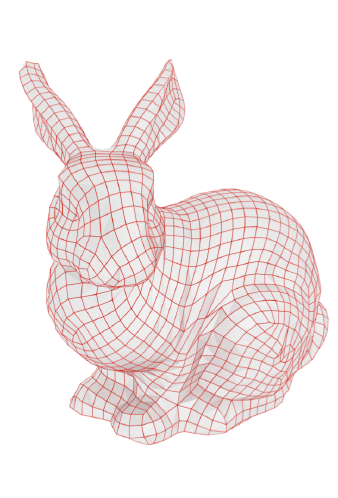}
			\label{fig:bunnyMesh}
		}
		\hfill
	\subfloat[Immersed lattice]
		{
			\includegraphics{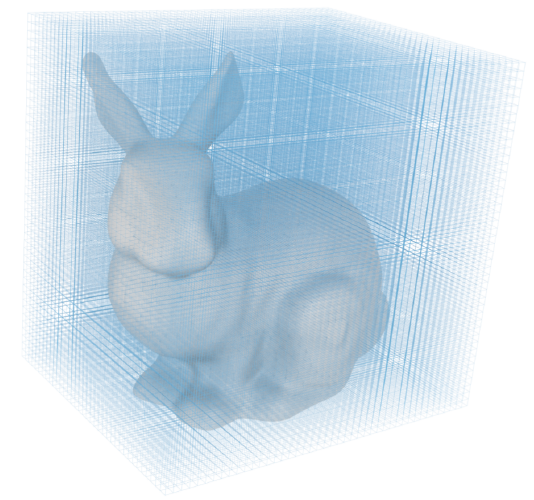}
		}
		\hfill
	\subfloat[Intersection points]
		{
			\includegraphics{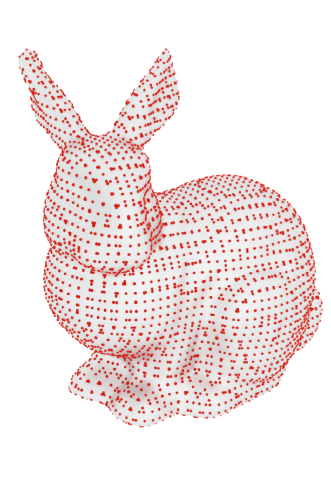}
			\label{fig:bunnyIntersect}
		}
\end{minipage}
\caption{Generation of an internal lattice for the Catmull--Clark subdivision surface of the Stanford bunny.}
\label{fig:bunny}
\end{figure}

%
In our first example we investigate the time needed for generating a lattice within a given surface. Figure~\ref{fig:bunny} shows the steps involved in  generating an internal lattice for the
Stanford bunny. The bunny surface is a Catmull--Clark subdivision surface and consists of $2071$ facets. The number of unit cells for the lattice is chosen with $44 \, \times44 \, \times36$. This leads to  $5704$ intersection points between the lattice and the surface. Each lattice cell is tessellated with struts to give a  body-centred cubic unit cell as  shown in Figure~\ref{fig:bunnyLattice}. The deformation of the lattice-skin coupled structure under a uniform loading is depicted in Figure~\ref{fig:bunnyDisp}.

\begin{figure}
\begin{minipage}{0.4\textwidth}
\centering
	\subfloat[Lattice orientation 0$^{\circ}$]
	{
		\includegraphics[scale=1.1]{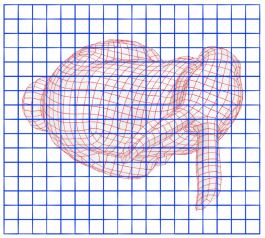}
		\label{fig:bunnyRotate0}
	}
	\vfill
	\subfloat[Lattice orientation 45$^{\circ}$]
	{
		\includegraphics[scale=1.1]{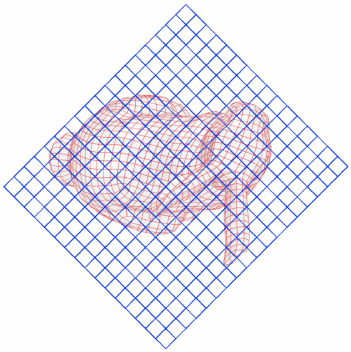}
		\label{fig:bunnyRotate45}
	}
\end{minipage}
\hfill
\begin{minipage}{0.6\textwidth}
\centering
	\subfloat[Computation times] 
	{
		\includegraphics[scale=0.45]{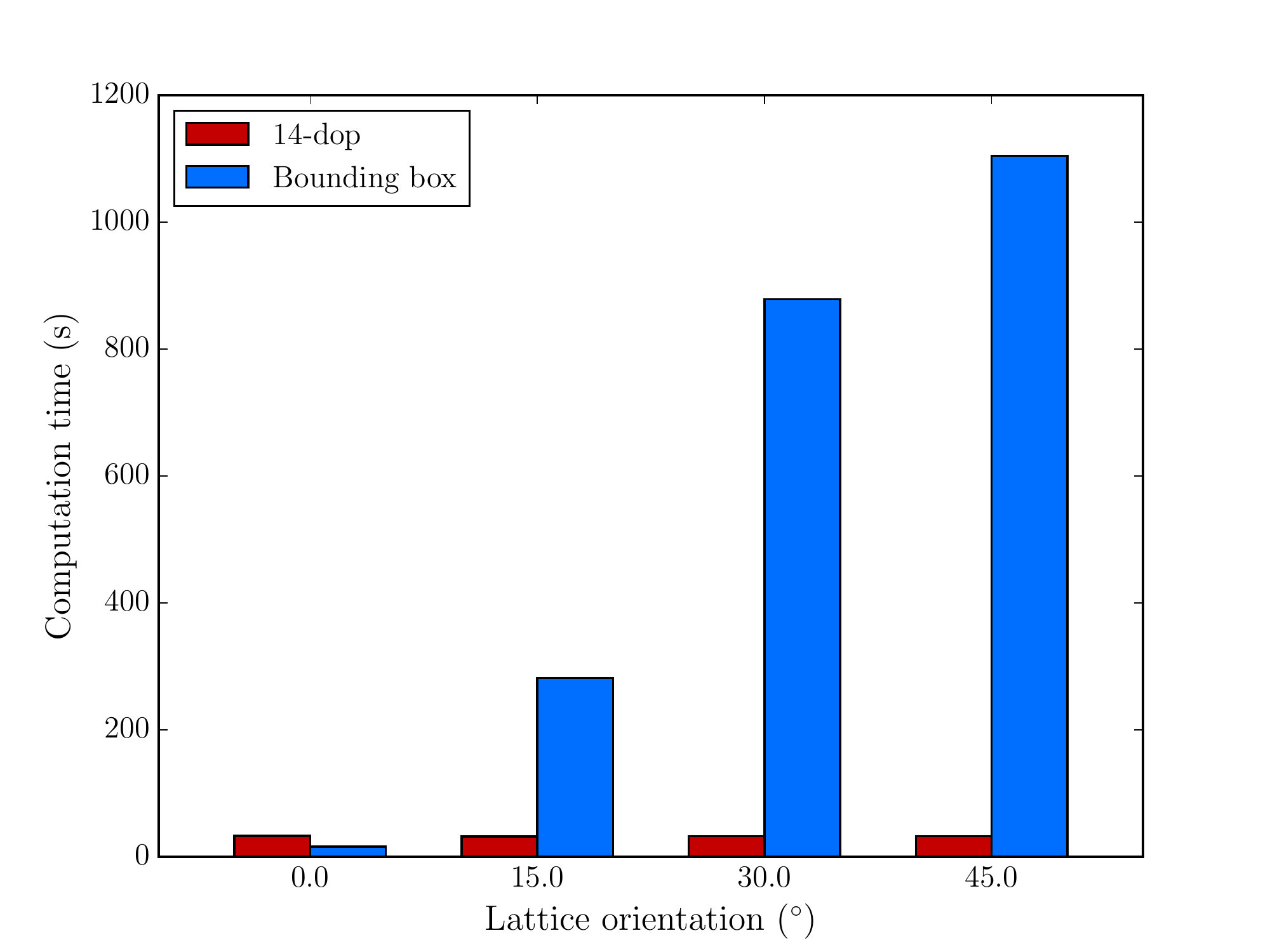}
	}
\end{minipage}
\caption{Efficiency comparison between k-dops and axis-aligned bounding boxes in dependence of the lattice orientation.   }
\label{fig:bunnyRotation}
\end{figure}

All the intersection points between the surface and lattice are computed using the intersection algorithm developed in Section~\ref{sec:geometry}. The SVD and eigenvalue computations involved are computed with the open source library Eigen~\cite{eigenweb}. For intersection detection we use 14-dops, consisting of six directions orthogonal to the coordinate planes, two directions parallel to the average normal of all the B\'ezier patches (evaluated by sampling), and six directions parallel to the lattice. There is no need to ensure that these 14 directions are always unique. In Figure~\ref{fig:bunnyRotation} the 14-dops and the axis-aligned bounding boxes (with $k = 6$) are compared for a lattice with $18 \, \times16 \, \times16$ unit cells in terms of efficiency of  the intersection detection. Different lattice orientations $\{ 0^\circ, \, 15^\circ, \, 30^\circ, \, 45^\circ\}$ with respect to the coordinate axis have been considered. As to be expected, 14-dops are far more efficient than the axis-aligned bounding boxes when the lattice is not aligned with the coordinate axis. Lattice lines not aligned with the coordinate axis have much larger axis-aligned bounding boxes. 

\begin{figure}
\centering
\includegraphics[width=0.7\linewidth]{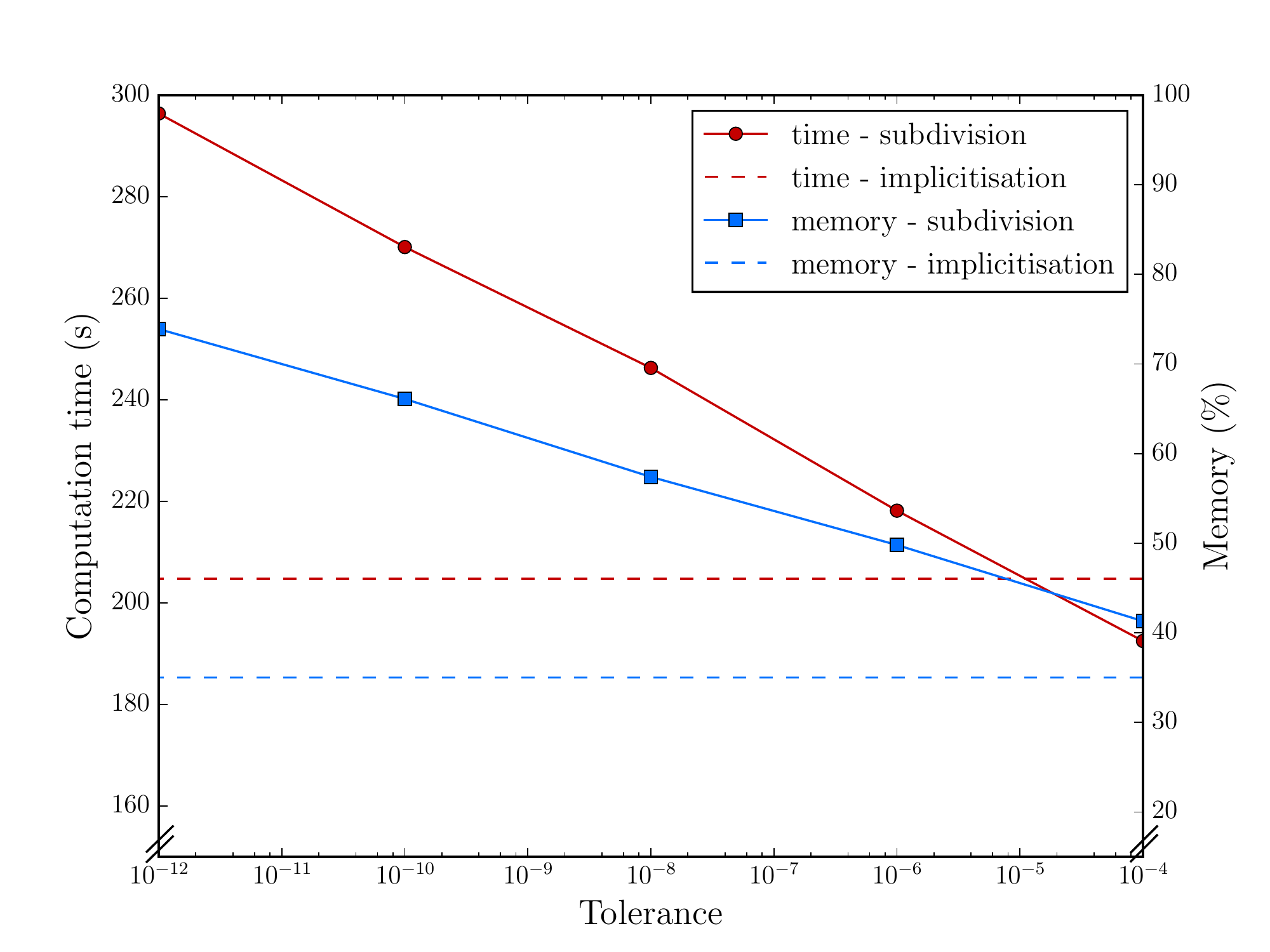}
\caption{Computation time and memory usage for the implicitisation and the subdivision methods for different flatness tolerances (lattice orientation is~$0^\circ$ and~14-dops are used).}
\label{fig:comparison}
\end{figure}

For evaluating the accuracy and efficiency of the  implicitisation, we compare it with the widely used subdivision, also referred to as divide-and-conquer, method for intersection computation \cite{houghton1985implementation}. In the subdivision method the B\'ezier patches are successively refined until they are sufficiently flat so that the intersection computation reduces to the trivial intersection between a triangle and a line segment.  In the chosen flatness criterion the difference between the maximum and minimum support heights $| h_{j,\max}^{\text{patch}} -  h_{j,\min}^{\text{patch}}|$ along the direction $j$ of the average normal of a B\'ezier patch are measured, cf.~\eqref{eq:supportHeight}. This difference controls the accuracy of the intersection computation and has to be less than a prescribed tolerance. The tolerance becomes, however, ineffective when the line is (almost) tangential to the surface, which is not a trivial problem for the subdivision method~\cite{hasle2007geometric}. The performance comparison between the implicitisation and the subdivision methods is shown in Figure~\ref{fig:comparison}. Both the computation time and memory consumed increase nearly linearly with the decreasing tolerance for the subdivision method. This reflects the linear relationship between the depth of the bounding volume tree and the required accuracy. In contrast, the implicitisation exhibits a high accuracy (more than $10^{-12}$) but takes less time and memory compared with the subdivision method. The computations were performed on a computer equipped with an Intel Xeon(R) CPU E5-2623 v3 @ 3.00GHz processor and 32GB memory. 
%

One advantage of the subdivision method is that it is, without any modification, applicable around extraordinary vertices, where each spline consists of a sequence of nested patches~\cite{Stam:1998aa,Peters:2008aa,zhang2018subdivision}. Therefore, in the present paper around extraordinary vertices the intersections are computed with the subdivision method. In the considered Catmull--Clark mesh for the Stanford bunny the subdivision method was applied in $266$ out of $2071$ patches always.

%
\subsection{Sandwich plate}
\begin{figure}
\begin{minipage}{0.5\textwidth}
\subfloat[Sandwich plate with pyramidal lattice core]
{
	\includegraphics[scale=1]{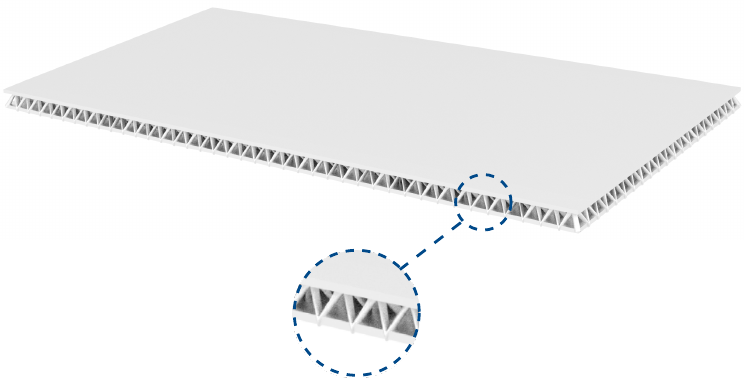}
	\label{fig:sandwichPlate}
}
\end{minipage}
\hspace{1cm}
\begin{minipage}{0.15\textwidth}
\subfloat[Unit cell with pyramidal tessellation]
{
	\includegraphics[scale=0.7]{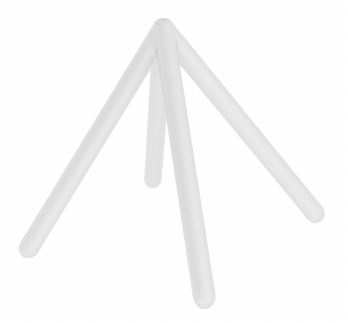}
	\hspace{0.5cm}
	\includegraphics[scale=0.7]{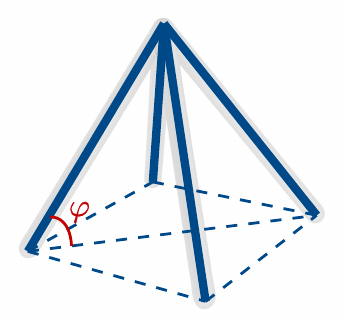}
	\label{fig:pyramdial}
}
\end{minipage}
\caption{Geometry of the sandwich plate.}
\end{figure}

%
Next we consider a sandwich plate with a pyramidal core to verify the accuracy and convergence of the finite element computations, see Figure~\ref{fig:sandwichPlate}. Although it is trivial to determine the intersection points between the flat skin plates and the lattice, the parametric coordinates of the intersection points still require the solution of nonlinear systems of algebraic equations. The size of the sandwich plate is chosen with~$1\, \mathrm{m} \, \times \, 2 \, \mathrm{m}$ and the distance between the top and bottom skins is~$0.05\, \mathrm{m}$. The pyramidal lattice core consists of cubic unit cells of side length~$0.05 \, \mathrm{m}$ and the struts are inclined by~$\varphi = \arctan(\!\sqrt{2})$ with respect to the bottom skin,  see Figure \ref{fig:pyramdial}.   The both skins and the struts are made of the same material with a Young's modulus~$E = 70\, \mathrm{GPa}$ and a Poisson's ratio $\nu=0.35$. The top skin is subjected to a uniform pressure loading of~$7\cdot10^5\, \mathrm{N}/\mathrm{m^2}$. The bottom skin is at its edges simply supported and each skin is discretised with~$256$ Kirchhoff--Love shell elements and cubic B-splines.

The homogenised material properties for the pyramidal lattice core can be found in~\cite{deshpande2001collapse} . Its in-plane Young's modulus is zero and its out-of-plane shear modulus is 
\begin{equation}
\overline{G} = \frac{\overline{\rho}}{8}E\sin^2 2\varphi \,  .
\end{equation}
The relative density~$\overline{\rho}$ of the core is 
\begin{equation}
	\overline{\rho} = \frac{\pi}{2 \cos^2\varphi\sin\varphi} \left (\frac{d}{l} \right)^2 \, ,
\end{equation}
where $d$ is the strut diameter and $l$ is the strut length.  Note that the two out-of-plane shear moduli have the same value~$\overline G$ due to the symmetry of the pyramidal core.  Analytic expressions for the displacements of a simply supported isotropic sandwich plate under uniform loading can be found in~\cite{allen1969analysis}. 

\begin{figure}
\centering
\subfloat[Comparison with analytic results]
{
	\includegraphics[scale=0.5]{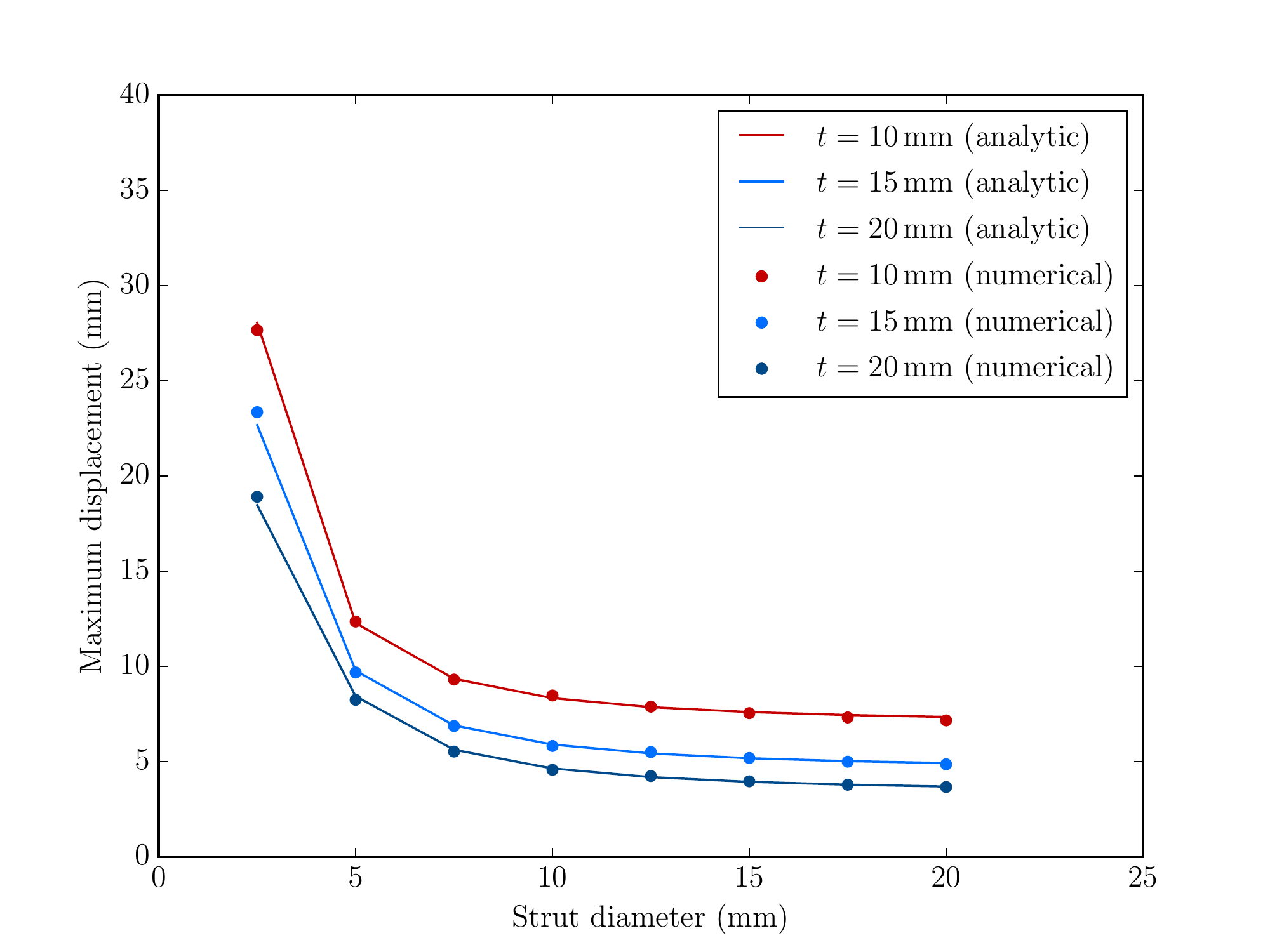}
	\label{fig:sandwichComparison}
}
\subfloat[Convergence of the maximum displacement  ($d = 17.5 \,  \mathrm{mm}$ and $t=15  \, \mathrm{mm}$) ]
{
	\includegraphics[scale=0.5]{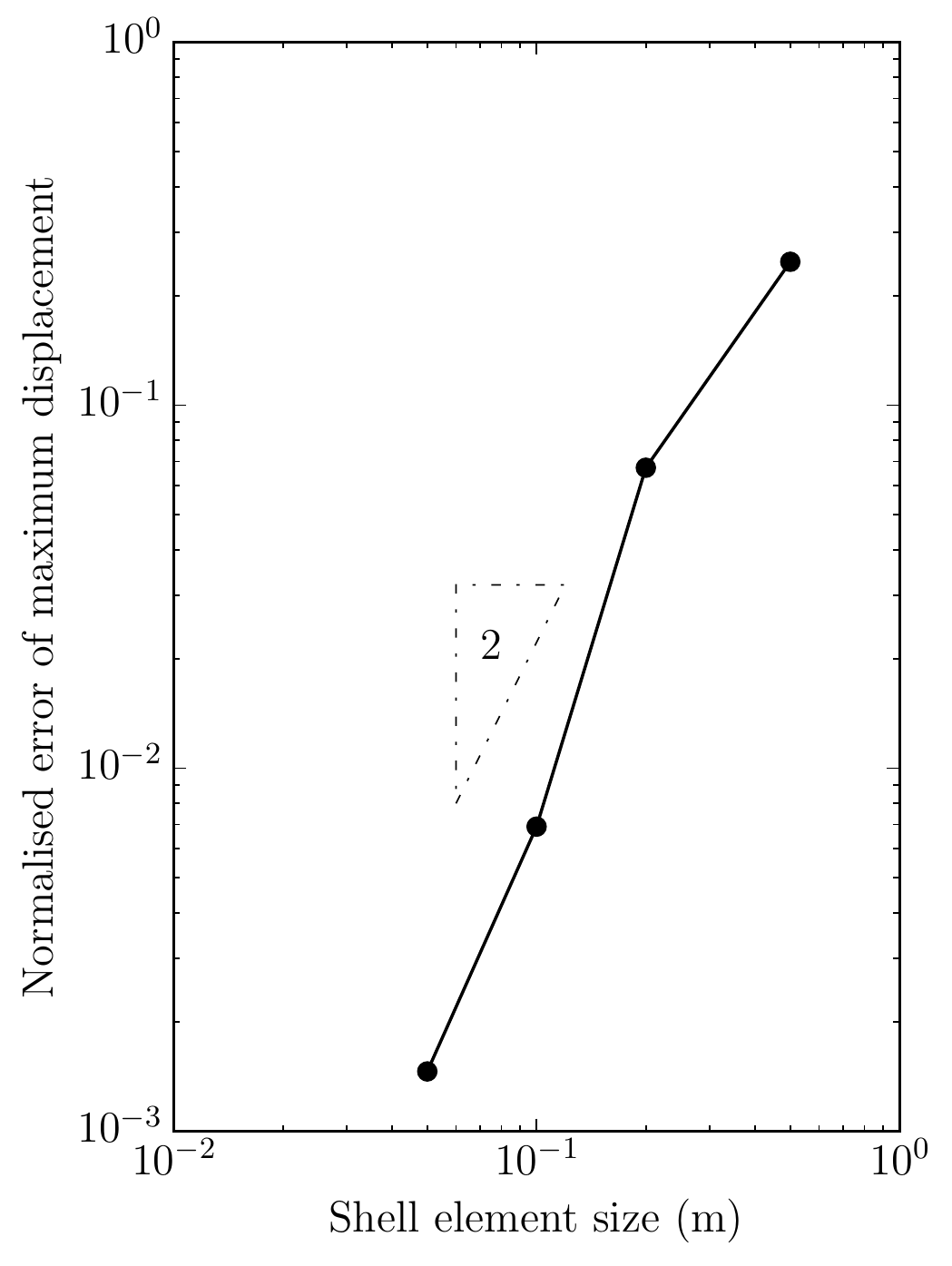}
	\label{fig:sandwichConvergence}
}
\caption{Maximum displacement of the sandwich plate under uniform pressure loading.}
\end{figure}

Three different skin thicknesses $t$ and eight different strut diameters $d$ are considered. The comparison of the numerical and the analytic maximum displacements is plotted in Figure~\ref{fig:sandwichComparison}. It can be seen that the numerical results agree well with the analytic ones for a range of different strut diameters and skin thicknesses. As depicted in Figure~\ref{fig:sandwichConvergence} the maximum displacement exhibits a quadratic convergence rate (for a plate with strut diameter $d = 17.5\,\mathrm{mm}$ and skin thickness $t= 15\,\mathrm{mm}$). The convergence rates for different strut diameters and skin thicknesses may fluctuate around the quadratic rate since in the reference solution~\cite{allen1969analysis} the lattice core is approximated with an isotropic bulk material, and the skin-lattice coupling has also an effect on the convergence rate.

\subsection{Spherical sandwich cap}
\begin{figure}
\centering
\includegraphics[scale=0.9]{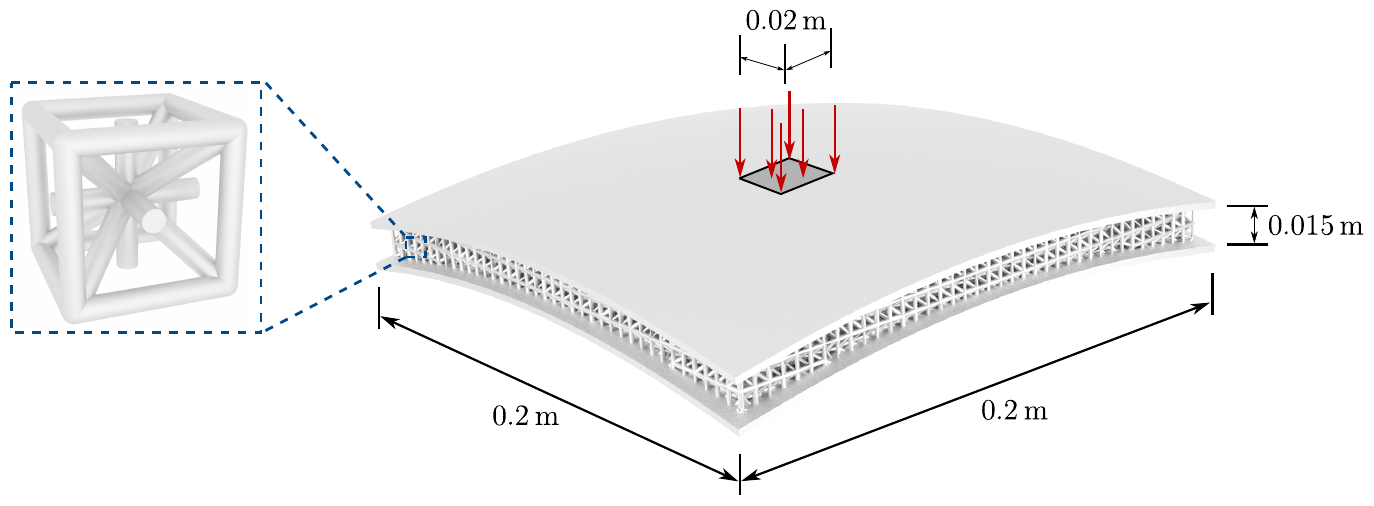}
\caption{Doubly-curved sandwich panel.}
\label{fig:curvedSandwich}
\end{figure}

%
As a further application of the developed analysis technique the optimal design of a spherical sandwich cap is considered. The sandwich cap has a BCC lattice core as shown in Figure~\ref{fig:curvedSandwich} and the lattice edges are aligned with the principal axes of the cap. Each unit cell has a side length of~$0.005\,\mathrm{m}$ and all struts have a diameter depending on the prescribed lattice volume ratio.  The size of the sandwich cap projected onto the horizontal plane is \mbox{$0.2\,\mathrm{m} \, \times \, 0.2\,\mathrm{m}$.} The distance between the top and bottom skin is $0.015\,\mathrm{m}$. Each skin is discretised with~$64$ Kirchhoff--Love shell elements and cubic B-splines. The Young's modulus and the Poisson's ratio of the thin-shell and the lattice material are $E = 100\,\mathrm{GPa}$ and $\nu=0$. A uniform loading with the magnitude $10^6\, \mathrm{N/m^2}$ is applied in the region $[-0.01\,\mathrm{m},  \, 0.01\,\mathrm{m}]\, \times \, [-0.01\,\mathrm{m}, \, 0.01\,\mathrm{m}]$ around the centre  of the upper skin. The four corners of the top and bottom skins are fixed.

\begin{figure}
\centering
\includegraphics[scale=0.55]{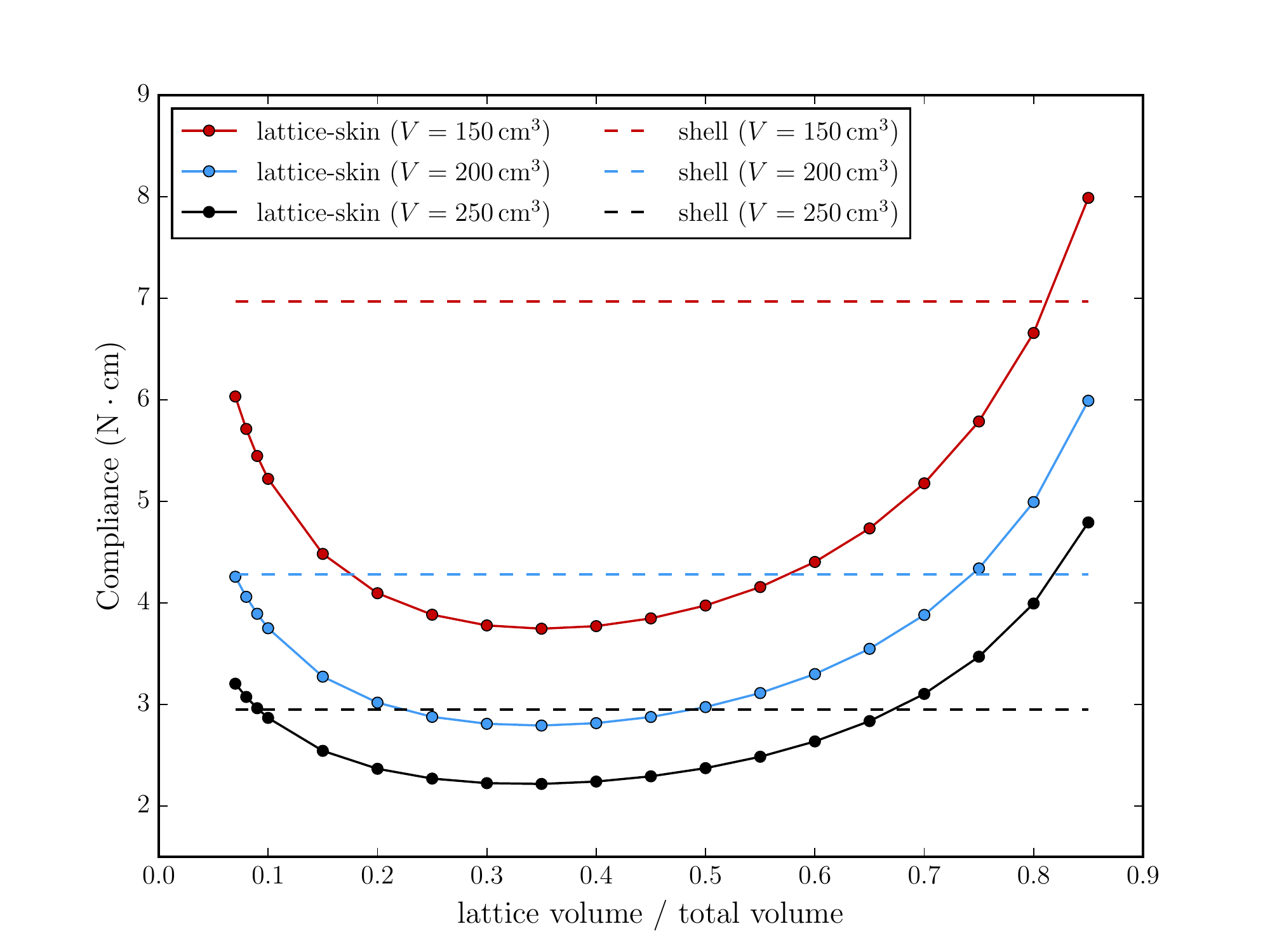}
\caption{Compliance for different lattice to total volume ratios.}
\label{fig:curvedSandwichComp}
\end{figure}

The total volume of the sandwich cap is prescribed and the optimal lattice to total volume ratio is sought. In Figure~\ref{fig:curvedSandwichComp} the structural compliance for different lattice to total volume ratios are plotted. Three different total volumes \mbox{$V \in \{ 150\, \mathrm{cm^3}, \, 200\, \mathrm{cm^3}, \, 250\, \mathrm{cm^3} \}$} are considered.  For each total volume there  exists an optimal lattice to total volume ratio as can be inferred from Figure~\ref{fig:curvedSandwichComp}. Structures consisting either only of a lattice or only one shell skin are non-optimal.  For a  relatively high lattice volume ratio the two shell skins become ineffective and the compliance becomes very large. In contrast, for a relatively low lattice  volume ratio the lattice is not sufficiently stiff to carry the shear forces occurring between the top and bottom skins. As also recently discussed in Sigmund et al.~\cite{sigmund2016non} and illustrated in Figure~\ref{fig:curvedSandwichComp} a lattice structure is not per se the most efficient structure. 

\section{Conclusions}
%

A surface interrogation technique for computing the intersection points of curves with spline surfaces has been introduced. The two key components are the hierarchical k-dop bounding volume tree for fast proximity search and the implicit matrix representation for intersection computations. As demonstrated, in contrast to conventional axis-aligned bounding boxes the k-dop bounding volumes can significantly shorten the search times and are, most crucially, insensitive to the orientation of the considered geometry. The matrix-based implicitisation of surface patches reduces intersection computations to the solution of SVD and generalised eigenvalue problems. The obtained non-complex eigenvalues represent all the intersection points between a surface patch and a curve. Multiple real eigenvalues indicate multiple intersections between the curve and the patch. While establishing the implicitisation matrix and solving the generalised eigenvalue problem the null spaces of a range of matrices need to be computed. The largest matrix is of size \mbox{$6 p^2 \, \times \, 8 p^2$,} \mbox{i.e. $54 \, \times \, 72$} for the cubic patches used in this paper. All the required null spaces and eigenvalues are computed with a high-performance linear algebra software library. The main advantage of the proposed interrogation technique is its robustness in comparison to Newton--Raphson iteration, or its variants, and its efficiency in comparison to subdivision based methods.  In none of the introduced examples, with up to $\approx 5500$ intersections, were any robustness issues observed. Subdivision based methods are not competitive due to strict accuracy requirements in finite element computations. This leads to very deep bounding volume trees, with $\gtrapprox 30$ levels, which may consume too much memory and take too long to construct and traverse. However, subdivision based methods can be applied without modification to almost any surface representation, which is relevant for interrogating close to extraordinary vertices. For instance, in the considered Catmull-Clark subdivision scheme each of the facets attached to an extraordinary vertex consists of a sequence of nested B\'ezier patches.  Although in each of the B\'ezier patches a matrix-based interrogation technique can be applied, the intrinsic need for subdivision refinement makes a subdivision based interrogation more appealing.

The developed interrogation technique has been applied to the isogeometric design and analysis of lattice-skin structures to be manufactured with 3d printing. The isogeometric design and analysis of such structures is currently hampered by the poor scalability of prevalent CAD systems and their limitations in representing non-solid entities, like lattices. The implemented approach combines a lattice, consisting of a set of nodes and their connectivity, with a surface for the skin. The lattice is discretised as a pin-jointed truss structure and the skin as a Kirchhoff--Love shell. This makes it feasible to design and analyse structures with several millions of entities and the entire approach is amenable to parallelisation. Although not pursued in this paper, the presented interrogation technique is also suitable for slicing the geometry model to generate machine instructions for the 3d printer~\cite{ma2004nurbs, starly2005direct}. 

In closing, we reiterate the importance of surface interrogation in a gamut of isogeometric analysis applications offering opportunities for research, including in enforcement of non-penetration constraints in contact, tessellating cut-cells in immersed or embedded boundary methods and coupling trimmed shell patches. We note also that it is straightforward to extend the introduced interrogation technique to other polynomial surface representations,  including triangular B\'ezier patches and Lagrange finite elements~\cite{laurent2014implicit, xiao2019non}. 

%
\section{Appendix}
%

%
\subsection{Products of Bernstein polynomials \label{sec:productAppendix}}
%
According to the definition of the univariate Bernstein polynomials~\eqref{eq:bernstein}, the product of two univariate Bernstein polynomials \mbox{$B_i^\lambda (\theta^1)$} with \mbox{$i \in \{ 1,  \, \dotsc, \, \lambda + 1 \}$} and \mbox{$B_j^\mu (\theta^1)$} with  \mbox{$j \in \{ 1, \, \dotsc, \, \mu + 1 \}$} can be expressed with
\begin{equation}
 	B_i^\lambda (\theta^1) B_j^\mu (\theta^1) = \frac{\binom{\lambda}{i - 1}\binom{\mu}{j - 1} }{\binom{\lambda + \mu}{i + j -2}}B_{i + j - 1}^{\lambda + \mu }(\theta^1)  \, .
\end{equation}
Similarly, the product of two bivariate Bernstein polynomials \mbox{$B_{\vec i}^{\vec \lambda} (\vec \theta) = B_{i^1}^{\lambda^1} (\theta^1) B_{i^2}^{\lambda^2} (\theta^2)$} and \mbox{$ B_{\vec j}^{\vec \mu} (\vec \theta)= B_{j^1}^{\mu^1}(\theta^1)B_{j^2}^{\mu^2}(\theta^2)$} can be expressed with
\begin{equation}
	B_{\vec i}^{\vec \lambda} (\vec \theta) B_{\vec j}^{\vec \mu} (\vec \theta)  = \frac{\binom{\lambda^1}{i^1 - 1}\binom{\mu^1}{j^1 - 1}\binom{\lambda^2}{i^2 - 1}\binom{\mu^2}{j^2 - 1}}{\binom{\lambda^1 + \mu^1}{i^1 + j^1 - 2}\binom{\lambda^2 + \mu^2}{i^2 + j^2 - 2}}B_{i^1 + j^1 - 1}^{\lambda^1 + \mu^1}(\theta^1)B_{i^2 + j^2 - 1}^{\lambda^2 + \mu^2}(\theta^2) \, .
\end{equation}

%
\subsection{Intersection of two linear B\'ezier curves \label{sec:intersectAppendix}}
%
In order to fix ideas, consider the trivial case of computing the intersection of two linear curves. The given B\'ezier curve $\vec{f}(\theta)$ has the control points \mbox{$\vec{x}_1 = (0, \, -1, \, 0)^\trans$} and \mbox{$\vec{x}_2 = (1, \, 1, \, 0)^\trans$} and the associated weights  \mbox{$w_1=w_2=1$}. The parametric representation of the curve in the homogeneous coordinates reads 
\begin{equation} \label{eq:linearBezier}
	\vec{f}(\theta) =
	(1 - \theta)
	\begin{pmatrix*}[r]
	0 \\
	-1 \\
	0 \\
	1
	\end{pmatrix*}
	+ \theta
	\begin{pmatrix*}
	1 \\
	1 \\
	0 \\
	1
	\end{pmatrix*}
	=
	\begin{pmatrix*}[c]
	\theta \\ 2\theta - 1 \\ 0 \\ 1
	\end{pmatrix*} \, .
\end{equation}

According to Bus\'e et al.~\cite{laurent2014implicit}, for a linear B\'ezier curve~$\vec f (\theta)$ the smallest possible polynomial degree for the auxiliary vector $\vec{g}(\theta)$ is zero. Hence, first, in Section~\ref{sec:intersectConstAppendix}, the intersection is computed with a constant auxiliary vector $\vec{g}(\theta)$. As discussed in Section~\ref{sec:implicitisation}, choosing a high-order polynomial~$\vec {g}(\theta)$ is possible. Increasing the polynomial degree leads, however, to larger matrices and make the intersection computations less efficient. With the sole aim of demonstrating the process of the generalised eigenvalue computation, in  Section~\ref{sec:intersectLinAppendix} the intersection is computed with a linear auxiliary vector $\vec{g}(\theta)$. 

%
\subsubsection{Constant auxiliary vector \label{sec:intersectConstAppendix}}
%
The constant auxiliary vector for computing the implicit matrix representation of the linear B\'ezier curve~\eqref{eq:linearBezier} is given by
\begin{equation} \label{eq:constAux}
	\renewcommand\arraystretch{1.2}
	\vec{g}(\theta) = \vec{\gamma}_1 = 
	\begin{pmatrix*}
	\gamma_{1}^{1} \\
	\gamma_{1}^{2} \\
	\gamma_{1}^{3} \\
	\gamma_{1}^{4}
	\end{pmatrix*} \, .
\end{equation}
According to the orthogonality requirement  \eqref{eq:orthogonality}, i.e.  $\vec{f}(\theta)\cdot\vec{g}(\theta) = 0$,  the two functions have to satisfy
\begin{equation}
	(1 - \theta)(-\gamma_{1}^{2} + \gamma_{1}^{4})  + \theta(\gamma_{1}^{1} + \gamma_{1}^{2}+ \gamma_{1}^{4})  = 0 \, , 
\end{equation}
or expressed in matrix form
\begin{equation} \label{eq:constOrthogonal}
	\renewcommand\arraystretch{1.2}
	\begin{pmatrix}
	1 - \theta & \theta
	\end{pmatrix}
	\begin{pmatrix}
	0 & -1 & 0 & 1 \\
	1 & 1 & 0 & 1
	\end{pmatrix}
	\begin{pmatrix}
	\gamma_{1}^{1} \\
	\gamma_{1}^{2} \\ 
	\gamma_{1}^{3} \\
	\gamma_{1}^{4}
\end{pmatrix} 
= 0 \, .
\end{equation}
By inspection, the left null space of the $2\times4$ matrix, denoted with $\vec C$ in Section~\ref{sec:implicitisation}, is spanned by the two vectors 
\begin{equation}
\vec{\gamma}_1^{(1)} = 
\begin{pmatrix}
-2 & 1 & 0 & 1
\end{pmatrix}^\trans
\,\quad \text{and} \quad
\vec{\gamma}_1^{(2)} = 
\begin{pmatrix}
0 & 0 & 1 & 0
\end{pmatrix}^\trans \, .
\end{equation}
These two left null vectors can also be readily obtained with a singular value decomposition or by bringing the matrix in a reduced row echelon form, see e.g.~\cite{strangLinAlg}. Introducing the null vectors into~\eqref{eq:constAux} yields the two auxiliary vectors
\begin{equation}
	\vec{g}^{(1)}(\theta) =
	\begin{pmatrix}
	-2 & 1 & 0 & 1
	\end{pmatrix}^\trans
	\,\quad \text{and} \quad
	\vec{g}^{(2)}(\theta) = 
	\begin{pmatrix}
	0 & 0 & 1 & 0
	\end{pmatrix}^\trans \, ,
\end{equation}
which yields, according to~\eqref{eq:planesMatrix}, the following implicit matrix representation
\begin{equation} \label{eq:implicitMatrix0}
	\vec{M} (\vec x) =
	\begin{pmatrix}
	-2x^{1} + x^{2} + 1 & x^{3}
	\end{pmatrix}.
\end{equation}
%
%
As an example, consider the intersection of $\vec{f}(\theta)$ with a second line 
\begin{equation} \label{eq:lineEquation}
	\vec{r}(\xi) =	
	(1 - \xi)
	\begin{pmatrix}
	0 \\ 1 \\ 0
	\end{pmatrix}
	+ \xi
	\begin{pmatrix}
	1 \\ 0 \\ 0
	\end{pmatrix}
	=
	\begin{pmatrix*}[c]
	\xi \\ 1 - \xi \\ 0
	\end{pmatrix*} \, .
\end{equation}
Substituting this line equation into the implicit matrix representation~\eqref{eq:implicitMatrix0} gives
\begin{equation}
	\vec{M} (\xi ) =
	\begin{pmatrix}
	-3 \xi  + 2 & 0
	\end{pmatrix}
	=
	\begin{pmatrix}
	2   & 0
	\end{pmatrix}
	+ \xi 
	\begin{pmatrix}
	-3   & 0
	\end{pmatrix}	
\end{equation}
At the intersection point~$\xi^*$ the matrix $\vec M(\xi)$ must be rank deficient, which is trivially the case when~$\xi^*=2/3$. Inserting~$\xi^*$ in~\eqref{eq:lineEquation} gives the intersection point $\vec r( \xi^*)=(2/3, \, 1/3, \, 0)$.

%
\subsubsection{Linear auxiliary vector \label{sec:intersectLinAppendix}}
\begin{figure}
\centering
\includegraphics[width=0.6\linewidth]{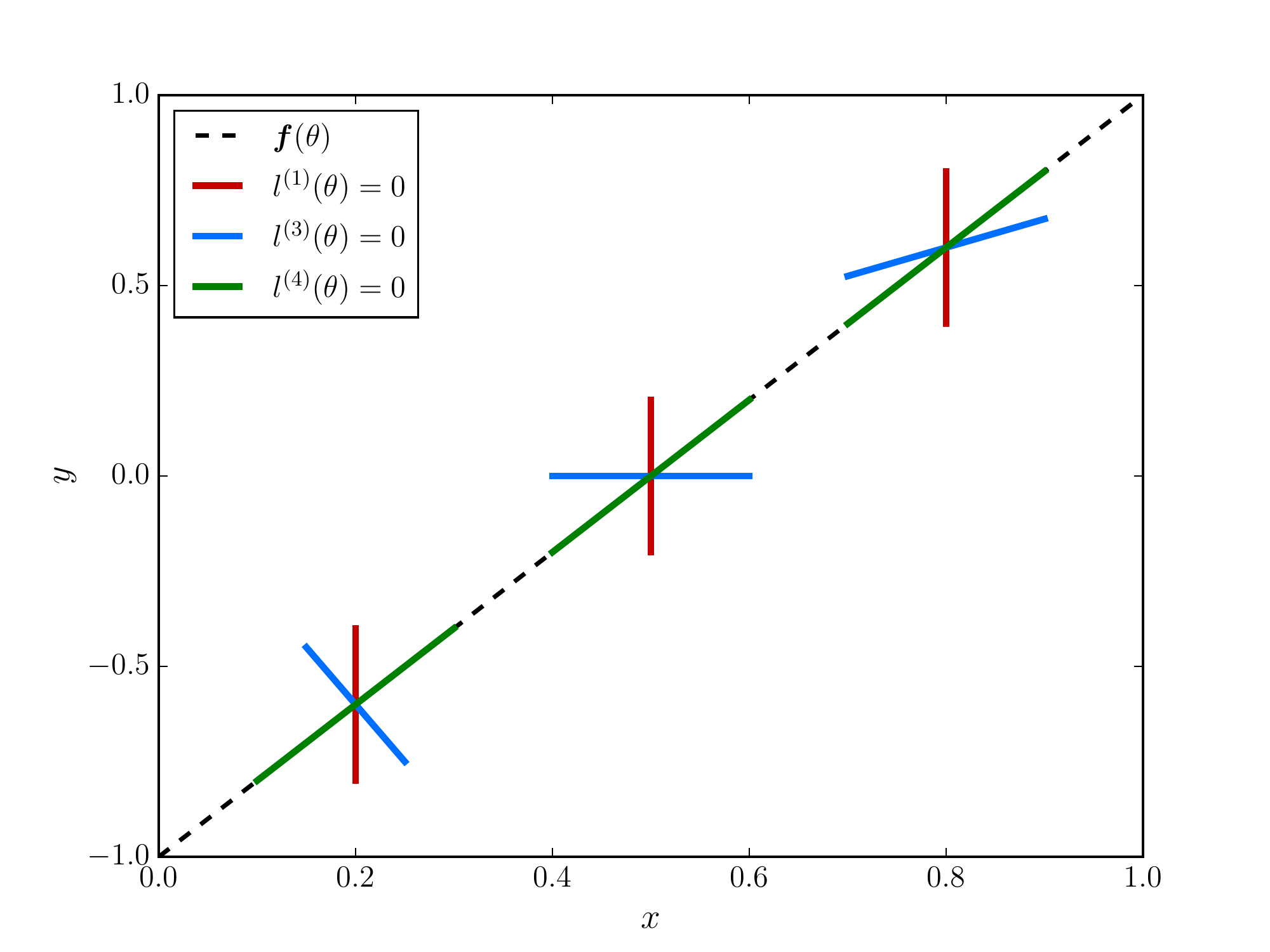}
\caption{The linear B\'ezier curve $\vec f(\theta)$ and the moving planes $l^{(i)}(\theta) = 0$ at $\theta \in \{0.2,\, 0.5,\, 0.8 \}$.  Note that the planes $l^{(2)}(\theta) = 0$ and $l^{(5)}(\theta) = 0$ are not shown because they correspond to the plane with $z = 0$.}
\label{fig:implicitLines}
\end{figure}

%
Next, the implicit matrix representation of the linear B\'ezier curve~\eqref{eq:linearBezier} is determined with the linear auxiliary vector
\begin{equation} \label{eq:auxiliaryBezier}
	\vec{g}(\theta) = (1 - \theta)\vec{\gamma}_1 + \theta\vec{\gamma}_2 \, ,
\end{equation}
where \mbox{$\vec{\gamma}_i = (\gamma_{i}^{1}, \, \gamma_{i}^{2}, \, \gamma_{i}^{3}, \, \gamma_{i}^{4})^\trans$}  with $i \in \{ 1, \, 2 \}$.  The orthogonality requirement \eqref{eq:orthogonality}, i.e.  \mbox{$\vec{f}(\theta)\cdot\vec{g}(\theta) = 0$,} yields
\begin{equation} \label{eq:orthogonal1}
	(1 - \theta)^2(-\gamma_1^2 + \gamma_1^4) + \theta(1 - \theta)(\gamma_1^1 +
	\gamma_1^2 - \gamma_2^2 + \gamma_1^4 + \gamma_2^4) + \theta^2(\gamma_2^1 +
	\gamma_2^2 + \gamma_2^4) = 0 \, ,
\end{equation}
or in matrix form
\begin{equation} \label{eq:orthogonalMatrix}
\begin{pmatrix}
(1 - \theta)^2 & 2\theta(1 - \theta) & \theta^2
\end{pmatrix}
\vec{C}
\vec{\gamma}
= 0 
\end{equation}
with
\begin{equation*} \label{coeffMatrix}
	\vec{C} = 
	\begin{pmatrix*}[r]
	0 & -1 & 0 & 1 & 0 & 0 & 0 & 0 \\
	\frac{1}{2} & \frac{1}{2} & 0 & \frac{1}{2} & 0 & -\frac{1}{2} & 0 & \frac{1}{2}
	\\
	0 & 0 & 0 & 0 & 1 & 1 & 0 & 1
	\end{pmatrix*} 
	 \qquad 
	 \text{and}
	 \quad
	 \vec \gamma = \begin{pmatrix}  \vec \gamma_{1} \\  \vec \gamma_{2} \end{pmatrix} \, .
\end{equation*}
It is worth pointing out that the single orthogonality requirement~\eqref{eq:orthogonality}, as can be inferred from~\eqref{eq:orthogonal1}, yields three equations for determining the eight components of the coefficients~$\vec{\gamma}_1$ and~$\vec{\gamma}_2$, which implies \mbox{$\dim (\ker \vec C) =5$}. The null space of $\vec C$, obtained by bringing it in row echelon form, is spanned by the five vectors
\begin{equation}  \label{eq:nullSpace}
\begin{aligned} 
	\vec{\gamma}^{(1)} &= 
	\begin{pmatrix}
	-1 & 0 & 0 & 0 & -1 & 0 & 0 & 1
	\end{pmatrix}^\trans \, , \qquad
	\vec{\gamma}^{(2)} &=
	\begin{pmatrix}
	\phantom{-} 0 & 0 & 0 & 0 & 0 & 0 & 1 & 0
	\end{pmatrix}^\trans \, , \\
	\vec{\gamma}^{(3)} &= 
	\begin{pmatrix}
	\phantom{-}1 & 0 & 0 & 0 & -1 & 1 & 0 & 0
	\end{pmatrix}^\trans \, , \qquad 
	\vec{\gamma}^{(4)} &= 
	\begin{pmatrix}
	-2 & 1 & 0 & 1 & 0 & 0 & 0 & 0
	\end{pmatrix}^\trans \, ,  \\
	\vec{\gamma}^{(5)} &= 
	\begin{pmatrix}
	\phantom{-}0 & 0 & 1 & 0 & \phantom{-}0 & 0 & 0 & 0
	\end{pmatrix}^\trans \, . 
	\end{aligned}
\end{equation}
The corresponding five auxiliary vectors follow from~\eqref{eq:auxiliaryBezier} with 
\begin{equation}
	\vec{g}^{(1)}(\theta) = 
	\begin{pmatrix*}[r]
	-1 \\ 0 \\  0 \\  \theta 
	\end{pmatrix*}
	\, , \quad
 	\vec{g}^{(2)}(\theta)  = 
	\begin{pmatrix*}[r]	
	0 \\ 0 \\ \theta \\ 0
	\end{pmatrix*}
	\, , \quad	
	\vec{g}^{(3)}(\theta) = 
	\begin{pmatrix*}[c]
	 1 - 2\theta \\  \theta \\ 0 \\  0
	 \end{pmatrix*}
	 \, , \quad
 	 \vec{g}^{(4)}(\theta)  =  
	 \begin{pmatrix*}[c]
	  2\theta - 2 \\  1-\theta \\ 0 \\  1-\theta
	  \end{pmatrix*}
	 \, , \quad
	\vec{g}^{(5)}(\theta) = 
	\begin{pmatrix*}[r]
	0 \\  0 \\ 1-\theta \\ 0
	  \end{pmatrix*} \, .
\end{equation}
As illustrated in Figure~\ref{fig:implicitLines}, the five vectors~$\vec{g}^{(i)}(\theta) $ can be associated for a given $  \theta^*$ with five implicitly defined planes 
\begin{equation} \label{eq:planeEquation}
	l^{(i)} ( \theta^* , \, \vec x) = \vec{g}^{(i)}( \theta^*)\cdot
	\begin{pmatrix}
	\vec x \\ 1
	\end{pmatrix} = 0 \, \quad \text{with } i \in \{1, \, \dotsc, \, 5 \} \, .
\end{equation}
Only at their intersection point~$\vec x^*$ all the plane equations $l^{(i)} ( \theta^*, \, \vec x^* ) = 0$ are simultaneously satisfied.  Due to the orthogonality condition \mbox{$\vec{f}(\theta^*)\cdot  \vec{g}^{(i)} (\theta^*) = 0$} the intersection point~$\vec x^*$ must also be the surface point
\begin{equation}
	\vec f ( \theta^* ) = \begin{pmatrix}  \vec x (\theta^*)   \\ 1   \end{pmatrix}  \, .
\end{equation}

The implicit matrix representation of~$\vec f(\theta)$ is obtained according to~\eqref{eq:planesMatrix} as follows 
\begin{align} \label{eq:planeMatrix}
\begin{pmatrix}
	\vec g^{(1)} (\theta) \cdot 
	 \begin{pmatrix} \vec x \\ 1 \end{pmatrix} 
	& \vec g^{(2)} (\theta)    \cdot 
	 \begin{pmatrix} \vec x \\ 1 \end{pmatrix} 
	& \dotsc
	 	& \vec g^{(5)} (\theta)    \cdot 
	 \begin{pmatrix} \vec x \\ 1 \end{pmatrix} 	
\end{pmatrix}
= 
\begin{pmatrix}
1 - \theta & \theta
\end{pmatrix}
\vec{M} \, 
\end{align}
with the implicitisation matrix 
\begin{equation} \label{eq:implicitMatrix}
	\vec{M} = 
	\begin{pmatrix}
	-x^{1} & 0 & x^{1} & -2x^{1} + x^{2} + 1 & x^{3} \\
	-x^{1} + 1 & x^{3} & -x^{1} + x^{2} & 0 & 0
\end{pmatrix}  \, .
\end{equation}
For a given $\theta^*$ and corresponding $\vec g (\theta^*)$ and $\vec f (\theta^*)$ the orthogonality condition \mbox{$\vec{f}(\theta^*)\cdot  \vec{g}^{(i)} (\theta^*) = 0$}  gives
\begin{equation} 
	\begin{pmatrix}
	1 - \theta^* &  \theta^*
	\end{pmatrix}
	\vec{M} (\theta^*)  = 
	\begin{pmatrix}
	0 & 0 & 0 & 0 & 0 
\end{pmatrix}
\,  .
\end{equation}
Hence,   $\vec M (\theta^*)$ is rank deficient and the vector $\begin{pmatrix} 1 - \theta^* &  \theta^* \end{pmatrix}$ lies in its  left null space. Indeed, at the point~$\vec{f}(\theta^*)$ the matrix has  $\rank (\vec M)=1$ and everywhere else it has $\rank (\vec M)=2$.  

As an example, consider now the intersection of $\vec{f}(\theta)$ with the line $\vec{r}(\xi)$ given in~\eqref{eq:lineEquation}.  Substituting~$\vec{r}(\xi)$ into the implicitisation matrix gives
\begin{equation}
	\vec{M}(\xi) = 
	\begin{pmatrix*}[c]
	-\xi & 0 & \xi & 2 - 3\xi & 0 \\	
	-\xi + 1 & 0 & 1 - 2\xi & 0 & 0
	\end{pmatrix*}
	=
	\vec{A} - \xi\vec{B} 
\end{equation}
with 
\begin{equation*}
	\vec{A} = 
	\begin{pmatrix*}[r]
	0 & 0 & 0 & 2 & 0 \\
	1 & 0 & 1 & 0 & 0
	\end{pmatrix*}
	\quad
	\text{and}
	\quad
	\vec{B} = 
	\begin{pmatrix*}[r]
	1 & 0 & -1 & 3 & 0 \\
	1 & 0 & 2 & 0 & 0
	\end{pmatrix*} \, .
\end{equation*}
The generalised eigenvalue problem~\eqref{eq:evKronecker}, i.e. $(\vec A - \xi \vec B)\vec \phi = \vec 0$,  for computing the intersection point~$\xi^*$ is solved by applying a sequence of orthogonal transformations as discussed in Section~\ref{sec:implicitisation}.  The first orthogonal transformation according to~\eqref{eq:columnCompression} gives the following matrices
\begin{equation}
	\vec{A}' = 
	\begin{pmatrix*}[r]
	-1.772680524 & 0.4369634904 & 0 & 0.8164965809 & 0 \\
	0.1438332239 & 1.346345178 & 0 & -0.4082482905 & 0
	\end{pmatrix*} \, ,
\end{equation}
\begin{equation}
	\vec{B}' =
	\begin{pmatrix*}[r]
	-3.297858703 & 0.3523180065 & 0 & 0 & 0 \\
	0.5351687938 & 2.171081381 & 0 & 0 & 0
\end{pmatrix*} \, .
\end{equation}
After applying two subsequent steps of orthogonal transformations $\vec{A}$ and $\vec{B}$ are reduced to the scalars
\begin{equation} \label{eq:reducedMatrix}
	\vec{A}''' = 
	\begin{pmatrix}
	1.549193338
	\end{pmatrix}
	\quad
	\text{and}
	\quad
	\vec{B}''' =
	\begin{pmatrix}
	2.323790008
	\end{pmatrix} \, .
\end{equation}
The generalised eigenvalue of matrices $\vec{A}'''$ and $\vec{B}'''$ is trivially \mbox{$\xi^* = 0.666666666$} and the coordinate of the intersection point is \mbox{$ \vec r (\xi^*) =\begin{pmatrix} \tfrac{2}{3} & \tfrac{1}{3} & 0 \end{pmatrix}^\trans$.}

At the intersection point, the implicitisation matrix becomes
\begin{equation}
	\renewcommand\arraystretch{1.2}
	\vec{M}(\xi^*) = 
	\begin{pmatrix*}[r]
	-\frac{2}{3} & 0 & \frac{2}{3} & 0 & 0 \\
	\frac{1}{3} & 0 & -\frac{1}{3} & 0 & 0
	\end{pmatrix*} \, 
\end{equation}
with its left null space spanned by the vector \mbox{$ \begin{pmatrix} \frac{1}{2} & 1 \end{pmatrix}$.} The parameter $\theta^*$ of the intersection point can be readily computed according to~\eqref{eq:planeMatrix} since the vector \mbox{$\begin{pmatrix} 1 - \theta^* &  \theta^* \end{pmatrix}$}  has to be collinear to \mbox{$\begin{pmatrix} \frac{1}{2} & 1 \end{pmatrix}$.}
Hence, $\theta^*= \frac{2}{3}$.

%
\subsection{Linearisation of implicit matrix representations \label{sec:linearisationAppendix}}
%
In Section~\ref{sec:implicitisation} only the intersection between a B\'ezier
patch and a line was considered. The introduced techniques can however easily be extended to the intersection between a B\'ezier patch and higher-order curves~\cite{thang2009curve}. As an example, consider the quadratic curve
\begin{equation}
	\vec r(\xi) = \vec c_2 \xi^2 + \vec c_1 \xi + \vec c_0
\end{equation}
with the prescribed vectors $\vec c_2$, $\vec c_1$ and $\vec c_0$. After substituting this curve, the implicit matrix representation~\eqref{eq:MRep} can be written as
\begin{equation}
	\vec{M}(\vec r( \xi)) = \vec M_0 + \vec M_1\xi + \vec M_{2}\xi^{2}  \, .
\end{equation}
As discussed, at the intersection points this matrix must be rank deficient, which motivates the generalised nonlinear eigenvalue problem 
\begin{equation}
	\vec{M}(\xi) \vec{\phi} = 
	 ( \vec M_0 + \vec M_1\xi + \vec M_{2}\xi^{2} ) \vec{\phi}
	= \vec 0 \, 
\end{equation}
which is equivalent to the eigenvalue problem 
\begin{equation}
	\left [ 
	\begin{pmatrix}
	\vec{0}& \vec{\vec I}  \\
	\vec{M}_0 & \vec{M}_1 
	\end{pmatrix} 
	- \xi 
	\begin{pmatrix}
	\vec{I} & \vec{0}  \\ 
	\vec{0} & -\vec{M}_2  \\
	\end{pmatrix}
	\right ] 
	\begin{pmatrix}
	\vec{\phi} \\  \xi \vec{\phi}
	\end{pmatrix}
	= \vec 0 \, .
\end{equation}
This generalised eigenvalue problem has the same form like $(\vec A - \xi \vec B)\vec \phi = \vec 0$ and can be solved in the same way, see Section~\ref{sec:implicitisation}. The obtained real eigenvalues  correspond to the intersection points of the quadratic curve with the patch.

\bibliographystyle{elsarticle-num-names}
\bibliography{latticeSkin}

\end{document}